\documentclass[11pt,a4paper,leqno]{article}
\usepackage{a4wide}
\usepackage{latexsym}
\usepackage{amsmath}
\usepackage{amssymb}
\newtheorem{defin}{Definition}
\newtheorem{lemma}{Lemma}
\newtheorem{prop}{Proposition}
\newtheorem{theo}{Theorem}
\newtheorem{corol}{Corollary}
\newenvironment{proof}{\medskip\par\noindent{\bf Proof}}{\hfill $\Box$
\medskip\par}
\begin{document}
\title{On $q-$asymptotics for linear $q-$difference-differential equations with Fuchsian and irregular singularities}
\author{Alberto Lastra, St\'ephane Malek, Javier Sanz}
\date{November, 14 2011}
\maketitle
\thispagestyle{empty}
{ \small \begin{center}
{\bf Abstract}
\end{center}
\noindent We consider a Cauchy problem for some family of linear $q-$difference-differential equations with Fuchsian and irregular singularities, that admit a unique formal power series solution in two variables $\hat{X}(t,z)$ for given formal power series initial conditions. Under suitable conditions and by the application of certain $q-$Borel and Laplace transforms (introduced by J.-P. Ramis and C. Zhang), we are able to deal with the small divisors phenomenon caused by the Fuchsian singularity, and to construct actual holomorphic solutions of the Cauchy problem whose $q-$asymptotic expansion in $t$, uniformly for $z$ in the compact sets of $\mathbb{C}$, is $\hat{X}(t,z)$. The small divisors's effect is an increase in the order of $q-$exponential growth and the appearance of a power of the factorial in the corresponding $q-$Gevrey bounds in the asymptotics.\medskip

\noindent Key words: $q-$difference-differential equations, $q-$Laplace transform, formal power series solutions, $q-$Gevrey asymptotic expansions, small divisors, Fuchsian and irregular singularities.\\
\noindent 2010 MSC: 34K25, 34M25, 34M30, 33E30.} \bigskip \bigskip

\section{Introduction}

The second author has studied partial differential equations of the form
\begin{equation}
t^{2r_{2}}\partial_{t}^{r_2}(z\partial_{z})^{r_{1}}\partial_{z}^{S}u(t,z) = F(t,z,\partial_{t},\partial_{z})u(t,z) \label{FIPDE}
\end{equation}
where $S,r_{1},r_{2}$ are nonnegative integers and $F$ is some differential operator with polynomial coefficients. These equations belong to a class of partial differential equations with both irregular singularity at $t=0$ (in the sense of T. Mandai~\cite{man2}) and Fuchsian singularity at $z=0$. These kind of problems have been extensively studied in the literature, see for example~\cite{al,bolepa,geta,ig,man,ya} for Fuchsian partial differential equations, and \cite{chluta,man2,ou1} for irregular singularities.

It is possible to construct formal power series solutions for the equation
(\ref{FIPDE}) of the form $\hat{u}(t,z) = \sum_{m \geq 0} \hat{u}_{m}(t) z^{m}/m!$, with coefficients in $\mathbb{C}[[t]]$, for given initial data
\begin{equation}
(\partial_{z}^{j}\hat{u})(t,0) = \hat{u}_{j}(t) \in \mathbb{C}[[t]], \label{FIPDE_init}
\end{equation}
$0\le j\le S-1$, which are assumed to be $1-$Borel summable with respect to $t$ in some direction $d \in \mathbb{R}$.

In the case $r_{1} = 0$, it was shown in \cite{ma4} that the formal series $\hat{u}(t,z)$ is $1-$Borel summable
with respect to $t$ in the direction $d$ if $d$ is well chosen, as series with coefficients in the Banach space of holomorphic functions near the origin with the supremum norm.

In the paper \cite{ma1} the case $r_{1} \neq 0$ was treated, noticing that the formal series $\hat{u}$ is in general no longer $1-$Borel summable, but one can show the existence of actual holomorphic solutions $u(t,z)$ which are Gevrey asymptotic of order larger than 1 to $\hat{u}(t,z)$ with respect to $t$ in sectors centered at 0 with finite radius in well chosen directions $d \in \mathbb{R}$. The reason for this different behaviour is the presence of small divisors introduced by the Fuchsian operator
$(z\partial_{z})^{r_{1}}$. More precisely, the singularities of the $1-$Borel transforms $\mathbb{B}\hat{u}_{m}(\tau)$ accumulate to the origin in $\mathbb{C}$ as $m$ tends to infinity, so that the $1-$Borel transform $(\mathbb{B}\hat{u})(\tau,z)$
with respect to $t$ is only holomorphic on a sector with infinite radius centered at $0$ with respect to $\tau$.
This induces a large multiplier effect on the $l-$th derivatives with respect to $t$ of the actual solution of
(\ref{FIPDE}) constructed with the classical Borel-Laplace procedure $u(t,z) = (\mathbb{L}\mathbb{B}\hat{u})(t,z)$ which grows like $CK^{l}l!\Gamma(1 + \gamma l)$ for some $C,K>0$ and some $\gamma > 2$ that can be expressed in terms of $r_{1}$,$r_{2}$, for all $l \geq 0$.
\medskip

In this paper, we study a $q-$analog of the problem (\ref{FIPDE}), (\ref{FIPDE_init}), discretized with respect to the variable
$t$, where $\partial_t$ is replaced by the operator $(f(qt)-f(t))/(qt-t)$ for $q\in\mathbb{C}$ (which formally tends to $\partial_t$ as $q$ tends to 1). Namely, we will consider the following linear $q-$difference-differential equation
\begin{equation}
( (z\partial_{z}+1)^{r_1}(t \sigma_{q})^{r_2} + 1)\partial_{z}^{S}\hat{X}(t,z) =
\sum_{k=0}^{S-1} b_{k}(z) (t\sigma_{q})^{m_{0,k}}(\partial_{z}^{k}\hat{X})(t, zq^{-m_{1,k}} ) \label{qFI}
\end{equation}
with given initial conditions
\begin{equation}
(\partial_{z}^{j}\hat{X})(t,0) = \hat{X}_{j}(t) \in \mathbb{C}[[t]],\qquad 0\le j\le S-1,
\label{qFI_init}
\end{equation}
where $S$,$m_{0,k}$,$m_{1,k}$ are nonnegative integers, for $0 \leq k \leq S-1$ and
where $q \in \mathbb{C}$ such that $|q|>1$, $\sigma_{q}$ is the dilation operator defined by
$(\sigma_{q}\hat{X})(t,z) = \hat{X}(qt,z)$, and $b_{k}(z)$ are polynomials in $z$.
As in previous works \cite{ma2}, \cite{ma3}, the map $(t,z) \mapsto (q^{m_{0,k}}t, zq^{-m_{1,k}})$ is assumed to be a volume shrinking map, meaning that the modulus of the Jacobian determinant $|q|^{m_{0,k} - m_{1,k}}$ is
less than~1. We will always assume that $r_2\ge 1$, while $r_1\ge 0$.

Advanced/delayed partial differential equations have also been widely studied, see for example ~\cite{ka,kaya1,kaya2,praspu1,yaka,zube}, and some authors have considered the use of special functions transforms for the study of the asymptotic properties of the solutions of $q-$difference-differential equations~\cite{Ho,praspu2}. Our present work is a contribution to this area.

It is not difficult to show (see Lemma~\ref{unique_formal_sol}) that this Cauchy problem has a unique formal power series solution of the form
$$
\hat{X}(t,z)=\displaystyle\sum_{h\ge 0}\hat{X}_h(t)\frac{z^h}{h!},
$$
where $\hat{X}_h(t)=\sum_{m\ge 0}f_{m,h}t^m\in\mathbb{C}[[t]]$, $h\ge 0$. Our purpose is to construct actual holomorphic solutions of this problem that are asymptotically represented by $\hat{X}(t,z)$ in a precise sense.

The key idea in our approach is the study of a related Cauchy problem,
\begin{equation}
((z\partial_{z} + 1)^{r_1} \tau^{r_2} + 1)
\partial_{z}^{S}\hat{W}(\tau,z) = \sum_{k=0}^{S-1} b_{k}(z)\tau^{m_{0,k}}(\partial_{z}^{k}\hat{W})(\tau,zq^{-m_{1,k}}) \label{CP2_intro}
\end{equation}
with initial conditions
\begin{equation}
(\partial_{z}^{j}\hat{W})(\tau,0) = \hat{W}_{j}(\tau) \in \mathbb{C}[[\tau]], \ \ 0 \leq j \leq S-1, \label{CP2_init_intro}
\end{equation}
which, by the application of a $q-$Laplace transform in the variable $\tau$, provides information on our initial problem (see Lemma~\ref{equiv_Cauchy_prob}). The $q-$Laplace transform we consider was introduced by J.-P. Ramis and C. Zhang in~\cite{razh}, and in recent years it has been used with great success in the study of the asymptotic properties of solutions of $q-$difference equations, see~\cite{dirasazh}, in much the same way as the classical Laplace-Borel transform has been applied to the asymptotic study of formal solutions to differential equations and singular perturbation problems in the complex domain (see the works of W. Balser~\cite{bals,bals2}, B. Malgrange~\cite{malg}, J.-P. Ramis~\cite{ra2} or O. Costin~\cite{co}).\\
\noindent This new Cauchy problem (\ref{CP2_intro}), (\ref{CP2_init_intro}) is studied in two respects.

Firstly, assuming the initial conditions $W_j$ are holomorphic and
have $q$-exponential growth (of order 2) in a set
$Vq^{\mathbb{Z}}=\{vq^h:v\in V,\ h\in\mathbb{Z}\}$, $V$ being a
well chosen bounded open set in $\mathbb{C}\setminus\{0\}$, and with some restriction on the argument of $q$, we
prove in Theorem~\ref{theo345} that there exists a unique solution
of (\ref{CP2_intro}), (\ref{CP2_init_intro}), of the form
\begin{equation}\label{formal_solution}
W(\tau,z) = \sum_{h \geq 0} W_{h}(\tau) \frac{z^h}{h!},
\end{equation}
holomorphic on $Vq^{\mathbb{Z}} \times \mathbb{C}$ and of $q$-exponential growth (of order 1) in $\tau$, in the terminology of~\cite{razh}, uniformly for $z$ in any compact set of $\mathbb{C}$. The increase in the order may be seen as an effect of the small divisors appearing in the problem.

Secondly, assuming the initial conditions $W_j$, $0\le j\le S-1$, are holomorphic near the origin, we prove in Theorem~\ref{teoHCP2} that the solution in (\ref{formal_solution}) has coefficients $W_h$ holomorphic in discs $D_h$ whose radii tend to 0 as $h$ tends to infinity, in such a way that there exist constants $C_1,T_1,X_1>0$ such that
$$\sup_{\tau\in \overline{D}_{j}}|\partial^{n}W_{j}(\tau)|\le
C_1\Big(\frac{1}{T_1}\Big)^{n}\Big(\frac{1}{X_1}\Big)^{j}n!j!(j+1)^{\frac{r_{1}n}{r_{2}}}|q|^{-j^2/2},
$$
for every $n,j\ge0$. The important fact here is the $q-$exponential decrease of these bounds with respect to $j$, what will turn out to be crucial in the following. These two results allow us to analyze the $q-$asymptotic expansion of the $q-$Laplace transforms of the $W_h$ (Proposition~\ref{q_Laplace_W_h}), which is shown to hold in a common domain $\mathcal{T}_{\lambda,q,\delta,r_{0}}$ (see~(\ref{disc_spiral}) for its definition) for all $h\ge 0$.

We are prepared to turn now to our main objective. Departing from formal initial conditions $\hat{X}_j$, $0\le j\le S-1$, whose $q-$Borel transforms $W_j(\tau)$ (in the terminology of~\cite{razh}) satisfy all the conditions in the previous two results, we are finally able to find a solution of our problem (Theorem~\ref{teo_solution}) in the form
$$
X(t,z)=\displaystyle\sum_{h\ge 0}\mathcal{L}_{q}^{\lambda}(W_h)(t)\frac{z^h}{h!},
$$
which is holomorphic in $\mathcal{T}_{\lambda,q,\delta,r_{0}}\times \mathbb{C}$, and such that given $R>0$, there exist constants $\tilde{C}>0$, $\tilde{D}>0$ such that for every $n\in\mathbb{N}$, $n\ge 1$, one has
$$
\Big|X(t,z)-\sum_{h\ge 0}\sum_{m=0}^{n-1}f_{m,h}t^m\frac{z^h}{h!}\Big|\le \tilde{C}\tilde{D}^n\Gamma(\frac{r_1}{r_2}(n+1))|q|^{n(n-1)/2}|t|^n
$$
for every $t\in\mathcal{T}_{\lambda,q,\delta,r_{0}}$, $z\in D(0,R)$. Again one may note that the small divisors phenomenon has caused the appearance of the term $\Gamma(\frac{r_1}{r_2}(n+1))$.

The paper is organized as follows. Section~\ref{sect_q_Laplace} provides the facts concerning the $q-$Laplace transform. Section~\ref{sectweightedBspaces} is devoted to the study of a first auxiliary Cauchy problem in suitable weighted Banach spaces of formal Laurent series. This is needed in the following section, devoted to the proof of Theorem~\ref{theo345}. A second Cauchy problem in weighted Banach spaces of formal Taylor series (Section~\ref{sect_second_CP}) is applied in the next Section, which contains Theorem~\ref{teoHCP2}. Finally, Section~\ref{sect_final_solution} consists of the construction of the solution, and it also contains some remarks on the nature of the solution in the special case that $r_1=0$, in which no small divisors appear.

We fix some conventions. $\mathbb{C}^{\ast}$ stands for $\mathbb{C}\setminus\{0\}$,
and $\mathbb{N}$ for the set $\{0,1,2,\cdots\}$. $D(0,r)$ denotes the open disc with center 0 and radius $r>0$. Given a set $V\subset\mathbb{C}$ and $q\in\mathbb{C}$, we define
$$
Vq^{\mathbb{Z}}=\{vq^h:v\in V,\ h\in\mathbb{Z}\},\quad Vq^{\mathbb{N}}=\{vq^h:v\in V,\ h\in\mathbb{N}\}.
$$

\section{A $q-$analogue of the Laplace transform and $q-$asymptotic expansion}\label{sect_q_Laplace}

In this section, we recall the definition of a $q-$analogue of the
Laplace transform introduced in the papers \cite{razh,zhang} and some of
its properties that will be useful in the sequel. For the sake of clarity, we include the proof of these results (mainly available in \cite{zhang}), since they contain important estimates that will be used in the proof of our main result (Theorem~\ref{teo_solution}).

\begin{prop}\label{proposition1}
Let $q \in \mathbb{C}$ such that $|q|>1$. Let $V$ be an open and bounded set in $\mathbb{C}^{\ast}$ and
$D(0,\rho_{0})$ a disc such that $V \cap D(0,\rho_{0}) \neq \emptyset$. Let $(\mathbb{F},|| . ||_{\mathbb{F}})$ be a
complex Banach space. Let $\phi : Vq^{\mathbb{N}} \cup D(0,\rho_{0}) \rightarrow \mathbb{F}$ be a holomorphic function
which satisfies the following estimates : there exist $C,M>0$ such that
\begin{equation}
||\phi(xq^{m})||_{\mathbb{F}} \leq M |q|^{m^{2}/2}C^{m} \label{|phi|<ql2}
\end{equation}
for all $m \geq 0$, all $x \in V$. Let $\Theta$ be the Jacobi Theta function defined in $\mathbb{C}^{\ast}$ by
$$ \Theta(x) = \sum_{n \in \mathbb{Z}} q^{-n(n-1)/2 } x^{n}. $$
Let $\delta >0$ and $\lambda \in V \cap D(0,\rho_{0})$. We denote by
\begin{equation}\label{disc_spiral}
\mathcal{R}_{\lambda,q,\delta} = \{ t \in \mathbb{C}^{\ast} : |1 + \frac{\lambda}{tq^{k}}| > \delta,
\forall k \in \mathbb{Z} \},    \quad  \mathcal{T}_{\lambda,q,\delta,r_{1}} = \mathcal{R}_{\lambda,q,\delta}\cap D(0,r_1).
\end{equation}
The $q-$Laplace transform of $\phi$ in the direction
$\lambda q^{\mathbb{Z}}$ is defined by
$$ \mathcal{L}_{q}^{\lambda}(\phi)(t) := \sum_{m \in \mathbb{Z}}
\phi(q^{m}\lambda)/ \Theta( \frac{q^{m}\lambda}{t})
$$
for all $t \in \mathcal{T}_{\lambda,q,\delta,r_{1}}$, if $r_{1} < |\lambda q^{1/2}|/C$. Moreover,
$\mathcal{L}_{q}^{\lambda}(\phi)(t)$ defines a bounded holomorphic function on $\mathcal{T}_{\lambda,q,\delta,r_{1}}$ with values in $\mathbb{F}$ when
$r_{1} < |\lambda q^{1/2}|/C$. Assume that the function $\phi$ has the following Taylor expansion
\begin{equation}
\phi(\tau) = \sum_{n \geq 0} \frac{ f_{n} }{ q^{n(n-1)/2} } \tau^{n} \label{taylor_phi}
\end{equation}
on $D(0,\rho_{0})$, where $f_{n} \in \mathbb{F}$, $n \geq 0$.
Then, there exist two constants $D,B>0$ such that
\begin{equation}
||\mathcal{L}_{q}^{\lambda}(\phi)(t) - \sum_{m = 0}^{n-1}f_{m}t^{m}||_{\mathbb{F}} \leq DB^{n}|q|^{n(n-1)/2}|t|^{n} \label{Laplace_asympt_expand}
\end{equation}
for all $n \geq 1$, for all $t \in \mathcal{T}_{\lambda,q,\delta,r_{1}}$.
\end{prop}

\noindent {\bf Remark:} In the situation described by~(\ref{Laplace_asympt_expand}) it is said that $\mathcal{L}_{q}^{\lambda}(\phi)$ admits the series $\sum_{m = 0}^{\infty}f_{m}t^{m}$ as $q-$Gevrey asymptotic expansion of order 1 (whenever the exponent of $|q|$ in the bounds is $n(n-1)/(2r)$ the order is said to be $r$). Analogously, a function that satisfies estimates such as (\ref{|phi|<ql2}) is said to have $q-$exponential growth of order 1 in $Vq^{\mathbb{N}}$.\\
\noindent If $\phi(z)=\sum_{n\ge 0}a_nz^n$ is an entire function such that there exists $C>0$ such that
$$
|a_n|\le C\exp(-(n-\alpha)^2/2)
$$
for all $n\ge 0$ and some $\alpha\ge 0$, then $\phi$ satisfies the estimates~(\ref{|phi|<ql2}). For a reference, see~\cite{ra}.\medskip

\begin{proof}
Since the Theta function $\Theta(x)$ satisfies the $q-$difference equation $\Theta(qx)=qx\Theta(x)$
for all $x \in \mathbb{C}^{\ast}$, we get that
\begin{equation}
\Theta( \frac{q^{m}\lambda}{t} ) = q^{m(m+1)/2} (\frac{\lambda}{t})^{m} \Theta( \frac{\lambda}{t} ) \label{prop1theta}
\end{equation}
for all $t \in \mathbb{C}^{\ast}$. Moreover, from Lemma 4.6 of \cite{rsz}, there exists
$K_{1}>0$ such that
\begin{equation}
 |\Theta(q^{m}\lambda/t)| \geq K_{1}\delta \sum_{n \in \mathbb{Z}} |q|^{-n(n-1)/2} |\frac{q^{m}\lambda}{t}|^{n}
\label{prop2theta}
\end{equation}
for all $t \in \mathcal{R}_{\lambda,q,\delta}$, all $m \in \mathbb{Z}$.\medskip

\noindent In the proof, we will show the estimates (\ref{Laplace_asympt_expand}). From them one may easily deduce that the series defining $\mathcal{L}_{q}^{\lambda}(\phi)(t)$ converges and defines a bounded holomorphic function on $\mathcal{T}_{\lambda,q,\delta,r_{1}}$. We would like to point out that many of the series following are initially formal, but we will finally prove their convergence.\\
Let $K \geq 0$ be an integer. First of all, we give estimates for the sum
$\sum_{m > 0} \phi(q^{m}\lambda)/\Theta(q^{m}\lambda/t)$. From the estimates (\ref{prop2theta}), we have that
\begin{equation}
|\Theta(\frac{\lambda}{t})| \geq K_{1}\delta |q|^{-K(K-1)/2}|\frac{\lambda}{t}|^{K} \label{|thetalambdat|>}
\end{equation}
for all $t \in \mathcal{R}_{\lambda,q,\delta}$. Using (\ref{|phi|<ql2}), (\ref{prop1theta}) and (\ref{|thetalambdat|>}), we get the estimates
$$
|| \frac{ \phi(q^{m}\lambda) }{ \Theta(q^{m}\lambda/t) } ||_{\mathbb{F}} \leq \frac{M}{K_{1}\delta}
(\frac{1}{|\lambda|})^{K} |q|^{K(K-1)/2}|t|^{K}(\frac{C|t|}{|\lambda||q|^{1/2}})^{m}
$$
for all $m >0$, all $t \in \mathcal{R}_{\lambda,q,\delta}$. So that if we choose a positive real number $r_{1} < |\lambda||q|^{1/2}/C$, then there exist
$D_{1},B_{1}>0$ (independent of $K$) such that
\begin{equation}
\sum_{m > 0} ||\frac{ \phi(q^{m}\lambda) }{ \Theta(q^{m}\lambda/t) }||_{\mathbb{F}} \leq
D_{1}(B_{1})^{K}|q|^{K(K-1)/2}|t|^{K} \label{|phi_m>0|<}
\end{equation}
for all $t \in \mathcal{T}_{\lambda,q,\delta,r_{1}}$.

\noindent In a second step, we give estimates for the sum $\sum_{m \leq 0} \phi(q^{m}\lambda)/\Theta(q^{m}\lambda/t) -
\sum_{n = 0}^{K} f_{n} t^{n}$, where the $f_{n}$ are defined in the Taylor expansion (\ref{taylor_phi}). From the formula
$$ q^{n(n-1)/2}t^{n} = \sum_{m \in \mathbb{Z}} \frac{ (q^{m}\lambda)^{n} }{ \Theta(q^{m}\lambda/t) } $$
for all $n \geq 0$, given in \cite{razh}, we can write (at least formally)
\begin{multline}
\sum_{m \leq 0} \phi(q^{m}\lambda)/\Theta(q^{m}\lambda/t) -
\sum_{n = 0}^{K} f_{n} t^{n} = \sum_{m \leq 0} \frac{1}{\Theta(q^{m}\lambda/t)}\left( \sum_{n \geq K+1} \frac{f_{n}}{q^{n(n-1)/2}}
(q^{m}\lambda)^{n} \right) \\
- \sum_{n=0}^{K} \frac{f_{n}}{q^{n(n-1)/2}} \left( \sum_{m > 0}
\frac{ (q^{m}\lambda)^{n} }{ \Theta(q^{m}\lambda/t) } \right) \label{decomp_sum_m<0}
\end{multline}
for all $t \in \mathbb{C}^{\ast}$. From the fact that $\phi$ has convergent expansion (\ref{taylor_phi})
on $D(0,\rho_{0})$, and since $|\lambda|<\rho_{0}$, there exist $C,A>0$, with $A<1/|\lambda|$, such that
\begin{equation}
||\frac{f_n}{q^{n(n-1)/2}}||_{\mathbb{F}} \leq CA^{n} \label{coeff_phi<}
\end{equation}
for all $n \geq 0$. From (\ref{decomp_sum_m<0}) and (\ref{coeff_phi<}), we deduce that
\begin{equation}
|| \sum_{m \leq 0} \phi(q^{m}\lambda)/\Theta(q^{m}\lambda/t) -
\sum_{n = 0}^{K} f_{n} t^{n} ||_{\mathbb{F}} \leq \mathcal{A}(t) + \mathcal{B}(t) \label{|phi_m<0 - polyn|<}
\end{equation}
where
$$ \mathcal{A}(t) = \sum_{m \leq 0}
\frac{1}{|\Theta(q^{m}\lambda/t)|}\left( \sum_{n \geq K+1} CA^{n}(|q|^{m}|\lambda|)^{n} \right) $$
and
$$ \mathcal{B}(t) = \sum_{n=0}^{K} CA^{n} \left( \sum_{m > 0}
\frac{ |(q^{m}\lambda)^{n}| }{ |\Theta(q^{m}\lambda/t)| } \right), $$
for all $t \in \mathbb{C}^{\ast}$.

\noindent  We give estimates for $\mathcal{A}(t)$. By summing up the geometric series (convergent because
$A|q|^m|\lambda|\le A|\lambda|<1$ for all $m\le 0$) and changing $m$ into $-m$, we first have that
there exists $D>0$ such that
\begin{equation}
\mathcal{A}(t) \leq DA^{K+1} \sum_{m \geq 0} \frac{ (|q|^{-m}|\lambda|)^{K+1} }{ |\Theta(q^{-m}\lambda/t)| } \label{A1<}
\end{equation}
for all $t \in \mathbb{C}^{\ast}$. From (\ref{prop2theta}), we have that
$$ |\Theta(q^{-m}\lambda/t)| \geq K_{1}\delta |q|^{-K(K-1)/2} |\frac{q^{-m}\lambda}{t}|^{K} $$
for all $m \geq 0$, all $t \in \mathcal{R}_{\lambda,q,\delta}$. We then have that
\begin{equation}
\frac{ (|q|^{-m}|\lambda|)^{K+1} }{ |\Theta(q^{-m}\lambda/t)| } \leq \frac{|\lambda|}{K_{1} \delta} |q|^{K(K-1)/2}|t|^{K}
(\frac{1}{|q|})^{m} \label{A2<}
\end{equation}
for all $m \geq 0$, all $t \in \mathcal{R}_{\lambda,q,\delta}$. From (\ref{A1<}) and (\ref{A2<}), we deduce that there
exist $D_{2},B_{2}>0$ (independent of $K$) such that
\begin{equation}
\mathcal{A}(t) \leq D_{2}(B_{2})^{K}|q|^{K(K-1)/2}|t|^{K} \label{calA<}
\end{equation}
for all $t \in \mathcal{R}_{\lambda,q,\delta}$.

\noindent In the next step, we get estimates for $\mathcal{B}(t)$.
From (\ref{prop2theta}), we have that
$$ |\Theta(q^{m}\lambda/t)| \geq K_{1}\delta |q|^{-(K+1)K/2} |\frac{q^{m}\lambda}{t}|^{K+1} $$
for all $m > 0$, all $t \in \mathcal{R}_{\lambda,q,\delta}$. We deduce that
\begin{equation}
\frac{ |(q^{m}\lambda)^{n}| }{ |\Theta(q^{m}\lambda/t)| } \leq \frac{ |\lambda|^{n} }{ K_{1}\delta }
(\frac{1}{|\lambda|})^{K+1}|q|^{(K+1)K/2}|t|^{K+1} (\frac{1}{|q|})^{m} \label{B1<}
\end{equation}
for all $m > 0$, all $0 \leq n \leq K$. From (\ref{B1<}), the equality $(K+1)K/2=K+K(K-1)/2$ and the fact that $|t|<r_1$ whenever $t \in \mathcal{T}_{\lambda,q,\delta,r_{1}}$, we obtain that there exist $D_{3},B_{3}>0$ (independent of~$K$) such that
\begin{equation}
\mathcal{B}(t) \leq D_{3}(B_{3})^{K}|q|^{K(K-1)/2}|t|^{K} \label{calB<}
\end{equation}
for all $t \in \mathcal{T}_{\lambda,q,\delta,r_1}$.\medskip

\noindent Finally, using the estimates
\begin{multline*}
 || \sum_{m \in \mathbb{Z}} \phi(q^{m}\lambda)/\Theta(q^{m}\lambda/t) -
\sum_{n = 0}^{K} f_{n} t^{n} ||_{\mathbb{F}} \leq || \sum_{m \leq 0} \phi(q^{m}\lambda)/\Theta(q^{m}\lambda/t) -
\sum_{n = 0}^{K} f_{n} t^{n} ||_{\mathbb{F}} \\
+ || \sum_{m > 0} \phi(q^{m}\lambda)/\Theta(q^{m}\lambda/t) ||_{\mathbb{F}}
\end{multline*}
we deduce from (\ref{|phi_m>0|<}), (\ref{|phi_m<0 - polyn|<}), (\ref{calA<}), (\ref{calB<}) that
\begin{equation}
|| \sum_{m \in \mathbb{Z}} \phi(q^{m}\lambda)/\Theta(q^{m}\lambda/t) -
\sum_{n = 0}^{K} f_{n} t^{n} ||_{\mathbb{F}} \leq D_4(B_4)^K|q|^{K(K-1)/2}|t|^{K}
\label{AmasB}
\end{equation}
for some $D_{4},B_{4}>0$ (independent of $K$). Now, for $K\in\mathbb{N}$, $K\ge 1$ one may write
$$
|| \sum_{m \in \mathbb{Z}} \phi(q^{m}\lambda)/\Theta(q^{m}\lambda/t) -
\sum_{n = 0}^{K-1} f_{n} t^{n} ||_{\mathbb{F}} \leq
|| \sum_{m \in \mathbb{Z}} \phi(q^{m}\lambda)/\Theta(q^{m}\lambda/t) -
\sum_{n = 0}^{K} f_{n} t^{n} ||_{\mathbb{F}} +|| f_{K} t^{K}||_{\mathbb{F}}
$$
and take into account (\ref{AmasB}) and (\ref{coeff_phi<}) in order to obtain
(\ref{Laplace_asympt_expand}), as desired.
\end{proof}

\begin{prop}\label{qLaplacetauphi} Let $V$ be an open and bounded set in $\mathbb{C}^{\ast}$ and $D(0,\rho_{0})$ be a disc such that
$V \cap D(0,\rho_{0}) \neq \emptyset$. Let $\phi$ be a holomorphic
function on $Vq^{\mathbb{N}} \cup D(0,\rho_{0})$ with values in $(\mathbb{F},||.||_{\mathbb{F}})$ which satisfies the estimates : There exist $C,K>0$, such that
\begin{equation}
||\phi(xq^{m})||_{\mathbb{F}} \leq K |q|^{m^{2}/2}C^{m} \label{|phi|<qm2Cm}
\end{equation}
for all $m \geq 0$, all $x \in V$. Then, the function
$M\phi(\tau) := \tau \phi(\tau)$ is holomorphic on $Vq^{\mathbb{N}} \cup D(0,\rho_{0})$ and satisfies estimates of the form (\ref{|phi|<ql2}). Let $\lambda \in V \cap D(0,\rho_{0})$. We have the following equality
$$
\mathcal{L}_{q}^{\lambda}(M\phi)(t) = t\mathcal{L}_{q}^{\lambda}(\phi)(qt)
$$
for all $t \in \mathcal{T}_{\lambda,q,\delta,r_{1}}$, if $r_{1} < |\lambda q^{1/2}|/(C|q|)$.
\end{prop}
\begin{proof} From the estimates (\ref{|phi|<qm2Cm}), we get a constant $r>0$ such that
$$ ||(M\phi)(xq^{m})||_{\mathbb{F}} \leq rK |q|^{m^{2}/2}(|q|C)^{m} $$
for all $m \geq 0$, all $x \in V$. From Proposition 1, we deduce that
$\mathcal{L}_{q}^{\lambda}(M\phi)(t)$ defines a holomorphic function on $\mathcal{T}_{\lambda,q,\delta,r_{1}}$, if $r_{1} < |\lambda q^{1/2}|/(C|q|)$. On the other hand,
\begin{equation}
t\mathcal{L}_{q}^{\lambda}(\phi)(qt) = \sum_{m \in \mathbb{Z}}
\frac{t\phi(q^{m}\lambda)}{\Theta(q^{m}\lambda/(qt))}. \label{def_tLq(phi)(qt)}
\end{equation}
But we have that
$$ \frac{t}{\Theta(q^{m}\lambda/(qt))} = \frac{q^{m}\lambda}{\Theta(q^{m}\lambda/t)} $$
for all $m \in \mathbb{Z}$. Indeed, put $y=q^{m}\lambda/(qt)$ in the identity $\Theta(qy) = qy\Theta(y)$.
From (\ref{def_tLq(phi)(qt)}), we get that
$$t\mathcal{L}_{q}^{\lambda}(\phi)(qt) = \sum_{m \in \mathbb{Z}}
\frac{q^{m}\lambda \phi(q^{m}\lambda)}{\Theta(q^{m}\lambda/t)} = \mathcal{L}_{q}^{\lambda}(M\phi)(t) $$
for all $t \in \mathcal{T}_{\lambda,q,\delta,r_{1}}$.
\end{proof}

For convenience, we recall the following concepts.
\begin{defin}\label{def_q_Borel}
A series $\hat{f}(t)=\sum_{n\ge 0}f_n t^n\in\mathbb{C}[[t]]$ is said to be $q-$Gevrey of order 1 if its so-called formal $q-$Borel transform of order 1,
$$
\hat{\mathcal{B}}_{q}\hat{f}(\tau)=\sum_{n \geq 0} \frac{ f_{n} }{ q^{n(n-1)/2} } \tau^{n},
$$
converges (i.e. it has positive radius of convergence).\\
The formal $q-$Laplace transform of order 1 of a series $\hat{g}(\tau)=\sum_{n\ge 0}g_n \tau^n\in\mathbb{C}[[\tau]]$ is defined as
$$
\hat{\mathcal{L}}_{q}\hat{g}(t)=\sum_{n \geq 0} q^{n(n-1)/2} g_{n} t^{n},
$$
so that these formal transforms are inverse of each other.
\end{defin}

It is immediate to check that, in agreement with Proposition~\ref{qLaplacetauphi}, we have that for every $\hat{g}\in \mathbb{C}[[\tau]]$,
\begin{equation}\label{property_q_Laplace_formal}
\hat{\mathcal{L}}_{q}(\tau\hat{g})(t)=t\hat{\mathcal{L}}_{q}\hat{g}(qt).
\end{equation}

\section{A Cauchy problem in a weighted Banach space of formal Laurent series}\label{sectweightedBspaces}

With the help of the $q-$Laplace transform we will change our initial problem (\ref{qFI}), (\ref{qFI_init}) into an equivalent one (\ref{CP2_intro}), (\ref{CP2_init_intro}), whose study will require the consideration of two auxiliary Cauchy problems. The first of them, which we are going to present in this Section, will be crucial in the study of the $q-$exponential growth of the coefficients of a solution of (\ref{CP2_intro}), (\ref{CP2_init_intro}). Although our equation involves a complex number $q$ with $|q|>1$, in this Section and in Section~\ref{sect_second_CP} we will be only concerned with the value $|q|$, so we directly work with a real value $q>1$.

\begin{defin}
We consider the vector space $\mathbb{E}_{q,(T,X)}$ of formal Laurent power series
\begin{equation}
 V(\xi,x) = \sum_{l \in \mathbb{Z},h \geq 0} v_{l,h} \xi^{l} \frac{x^h}{h!} \in \mathbb{C}[[\xi,\xi^{-1},x]] \label{V}
\end{equation}
such that
$$||V(\xi,x)||_{(T,X)} := \sum_{l \in \mathbb{Z},h \geq 0} \frac{|v_{l,h}|}{q^{P(l,h)}} T^{l} \frac{X^h}{h!} <\infty,$$
where $T,X>0$, $q > 1$ are positive real numbers and where
$$ P(l,h) = \left\{ \begin{array}{ll}
                     \frac{1}{4}l^{2} + \frac{1}{2}lh - \frac{1}{2}h^2 & \mbox{if $l \geq 0$,$h \geq 0$,} \\
                      - (1/2)h^{2} & \mbox{if $l \leq 0$,$h \geq 0$.}
                     \end{array} \right. $$
The space $(\mathbb{E}_{q,(T,X)}, ||.||_{(T,X)})$ is a Banach space.
\end{defin}

\noindent {\bf Remark:} Notice that we have a continuous inclusion $(\mathbb{E}_{q,(T,X')},||.||_{(T,X')})  \hookrightarrow (\mathbb{E}_{q,(T,X)}, ||.||_{(T,X)})$
when $0 < X \leq X'$.\medskip

\noindent We consider the integration operator $\partial_{x}^{-1}$ defined on $\mathbb{C}[[\xi,\xi^{-1},x]]$ by
$$ \partial_{x}^{-1}(V(\xi,x)) := \sum_{l \in \mathbb{Z},h \geq 1} v_{l,h-1} \xi^{l} \frac{x^h}{h!} \in \mathbb{C}[[\xi,\xi^{-1},x]] $$

\begin{lemma}\label{lema1} Let $m_1,s,h_{1},h_{2} \geq 0$ be nonnegative integers. Let $T,X>0$. Assume that the inequalities hold
\begin{equation}
s+ h_{2} \geq 2h_{1} \ \ , \ \ m_{1} \geq s + h_{2}. \label{ineq_shm}
\end{equation}
Then, there exist $C>0$ (depending $q$, $s,h_{1},h_{2},m_{1}$) such that
\begin{equation}
|| x^{s}(\partial_{x}^{-h_{2}}V)(q^{h_1}\xi, \frac{x}{q^{m_1}}) ||_{(T,X)} \leq CX^{(s+h_2)}
||V(\xi,x)||_{(T,X)} \label{||xTdxV||<C||V||}
\end{equation}
for all $V(\xi,x) \in \mathbb{E}_{q,(T,X)}$.
\end{lemma}
\begin{proof} Let $V(\xi,x) \in \mathbb{C}[[\xi,\xi^{-1},x]]$ as in~(\ref{V}). We have that
$$
x^{s}(\partial_{x}^{-h_{2}}V)(q^{h_1}\xi, \frac{x}{q^{m_1}}) = \sum_{l \in \mathbb{Z}, h\geq h_{2}+s}
v_{l,h-(s+h_{2})} \frac{ q^{h_{1}l} h!}{ q^{m_{1}(h-s)} (h-s)! } \xi^{l} \frac{{x}^{h}}{h!}
$$
From the definition of the norm $||.||_{(T,X)}$, we get that
\begin{multline}
|| x^{s}(\partial_{x}^{-h_{2}}V)(q^{h_{1}}\xi, \frac{x}{q^{m_1}}) ||_{(T,X)} =
\sum_{l \in \mathbb{Z}, h\geq h_{2}+s} \frac{ |v_{l,h-(s+h_{2})}| }{ q^{P(l,h-(s+h_{2}))} } T^{l}
\frac{ X^{h-(s+h_{2})} }{ (h-(s+h_{2})) ! } \times \\
\left\{ \frac{1}{ q^{P(l,h) - P(l,h-(s+h_{2})) - h_{1}l + m_{1}(h-s)} } \frac{ (h-(s+h_{2}))! }{ (h-s)! }
X^{s+h_{2}} \right\}. \label{factor_||xsTdxV||}
\end{multline}
In the rest of the proof, we will show that there exists a constant $C>0$
(depending on $q$, $s,h_{1},h_{2},m_{1}$) such that
\begin{equation}
\frac{1}{ q^{P(l,h) - P(l,h-(s+h_{2})) - h_{1}l + m_{1}(h-s)} } \leq C \label{|q^P|<C}
\end{equation}
for all $l \in \mathbb{Z}$, $h \geq 0$. Indeed, if $l \geq 0$, then
\begin{multline*}
P(l,h) - P(l,h-(s+h_{2})) - h_{1}l + m_{1}(h-s) = l( \frac{s+h_{2}}{2} - h_{1} ) + h( m_{1} - (s+h_{2})) \\
 - m_{1}s + \frac{(s+h_{2})^{2}}{2}
\end{multline*}
for all $h \geq 0$. From the assumption (\ref{ineq_shm}), we deduce that the inequalities (\ref{|q^P|<C}) hold
for $l \geq 0$, $h \geq 0$.

\noindent If $l \leq 0$, then
$$ P(l,h) - P(l,h-(s+h_{2})) - h_{1}l + m_{1}(h-s) = l(-h_{1}) + h( m_{1} - (s+h_{2})) -m_{1}s +
\frac{(s+h_{2})^{2}}{2} $$
for all $h \geq 0$. From the assumption (\ref{ineq_shm}), we deduce that the inequalities (\ref{|q^P|<C}) hold
for $l \leq 0$, $h \geq 0$.

\noindent Finally, the inequality (\ref{||xTdxV||<C||V||}) follows from the expression
(\ref{factor_||xsTdxV||}) and the estimates (\ref{|q^P|<C}).
\end{proof}

\begin{lemma}\label{lema2} Let $s,h_{1} \geq 0$ and $T_{0},X_{0} > 0$. Then, there exists a constant
$C_{1} > 0$ (depending on $q,s,h_{1},T_{0},X_{0}$) such that for all $0 < X_{1} \leq X_{0}q^{-s}$ and for all
$T_{1}>0$ satisfying
\begin{equation}
 q^{-h_{1}}T_{0} \leq T_{1} \leq T_{0}q^{\frac{s}{2} - h_{1}} \label{ineq_T0T1}
\end{equation}
one has
\begin{equation}
|| x^{s}V(q^{h_{1}}\xi,x) ||_{(T_{1},X_{1})} \leq C_{1} ||V(\xi,x) ||_{(T_{0},X_{0})} \label{|xTVT1X1|<|VT0X0|}
\end{equation}
for all $V(\xi,x) \in \mathbb{E}_{q,(T_{0},X_{0})}$.
\end{lemma}
\begin{proof} Let $V(\xi,x) \in \mathbb{C}[[\xi,\xi^{-1},x]]$. From the definition of the norm $||.||_{(T,X)}$, one can write
\begin{multline} || x^{s}V(q^{h_1}\xi,x) ||_{(T_{1},X_{1})} = \sum_{l \in \mathbb{Z}, h \geq s}
\frac{ |v_{l,h-s}| }{ q^{P(l,h-s)} } T_{0}^{l} \frac{ X_{0}^{(h-s)} }{ (h-s)! } \times \\
\left\{ \frac{1}{ q^{ P(l,h) - P(l,h-s) - h_{1}l } } (\frac{T_{1}}{T_{0}})^{l} (\frac{X_{1}}{X_{0}})^{h}
X_{0}^{s} \right\}
\label{xTV_TX_T0X0}
\end{multline}
In the rest of the proof, we will show that there exists a constant $C_{1}>0$ (depending on $q,s,h_{1}$) such that for all $0 < X_{1} \leq X_{0}q^{-s}$ and all $T_{1}$ satisfying
(\ref{ineq_T0T1}), one has
\begin{equation}
\frac{1}{ q^{ P(l,h) - P(l,h-s) - h_{1}l } } (\frac{T_{1}}{T_{0}})^{l} (\frac{X_{1}}{X_{0}})^{h} \leq C_{1} \label{|q^P|TlXh<C}
\end{equation}
for all $l \in \mathbb{Z}$, all $h \geq 0$.\medskip

\noindent If $l \geq 0$, then
$$ P(l,h) - P(l,h-s) - h_{1}l = l( \frac{s}{2} - h_{1} ) + h(- s ) + \frac{s^2}{2} $$
for all $h \geq 0$. So, we get that (\ref{|q^P|TlXh<C}) holds for all $l \geq 0$, all $h \geq 0$.\medskip

\noindent If $l \leq 0$, then
$$ P(l,h) - P(l,h-s) -h_{1}l = l( -h_{1} ) + h( -s) + \frac{s^2}{2} $$
for all $h \geq 0$. Hence, we get that
(\ref{|q^P|TlXh<C}) holds for all $l \leq 0$, all $h \geq 0$.\medskip

\noindent Finally, the inequality (\ref{|xTVT1X1|<|VT0X0|}) follows from the expression (\ref{xTV_TX_T0X0}) and the estimates
(\ref{|q^P|TlXh<C}).
\end{proof}

\begin{lemma}\label{lema3} Let $h_{2} \geq 0$ and $T_{0},X_{0}>0$. Then, there exists a constant $C_{2}>0$ (depending on $q,h_{2},T_{0},X_{0}$) such that
for all $0 < X_{1} \leq X_{0}q^{-h_{2}}$ and for all $T_{1}>0$ satisfying
\begin{equation}
T_{0} \leq T_{1} \leq T_{0}q^{h_{2}/2}, \label{ineq_T0T1_h2}
\end{equation}
one has
\begin{equation}
||\partial_{x}^{-h_{2}}V(\xi,x)||_{(T_{1},X_{1})} \leq C_{2}||V(\xi,x)||_{(T_{0},X_{0})} \label{|dxV|T1X1<|V|T0X0}
\end{equation}
for all $V(\xi,x) \in \mathbb{E}_{q,(T_{0},X_{0})}$.
\end{lemma}
\begin{proof} Let $V(\xi,x) \in \mathbb{C}[[\xi,\xi^{-1},x]]$. From the definition of the norm $||.||_{(T,X)}$, one can write
\begin{multline}
||\partial_{x}^{-h_{2}}V(\xi,x)||_{(T_{1},X_{1})} = \sum_{l \in \mathbb{Z}, h \geq h_{2}}
\frac{ |v_{l,h-h_{2}}| }{ q^{P(l,h-h_{2})} } T_{0}^{l} \frac{ X_{0}^{(h-h_{2})} }{ (h-h_{2})! } \times \\
\left\{ \frac{ 1 }{ q^{P(l,h) - P(l,h-h_{2})} } (\frac{T_{1}}{T_{0}})^{l} (\frac{X_{1}}{X_{0}})^{h}
\frac{ (h-h_{2})! }{ h! } X_{0}^{h_2} \right\} \label{dxV_T1X1_T0X0}
\end{multline}
In the rest of the proof, we will show that there exists a constant $C_{2}>0$ (depending on $q,h_{2}$) such that for all $0 < X_{1} \leq X_{0}q^{-h_{2}}$ and $T_{1}$ satisfying
(\ref{ineq_T0T1_h2}) one has
\begin{equation}
\frac{1}{ q^{ P(l,h) - P(l,h-h_{2}) } } (\frac{T_{1}}{T_{0}})^{l} (\frac{X_{1}}{X_{0}})^{h} \leq C_{2}
\label{|q^P|_h2TlXh<C}
\end{equation}
for all $l \in \mathbb{Z}$, all $h \geq 0$.\medskip

\noindent If $l \geq 0$, then
$$ P(l,h) - P(l,h-h_{2}) = l(\frac{h_{2}}{2}) + h(-h_{2}) + \frac{h_{2}^{2}}{2},$$
for all $h \geq 0$. Hence, we get that (\ref{|q^P|_h2TlXh<C}) holds for all $l \geq 0$, all $h \geq 0$.\medskip

\noindent If $l \leq 0$, then
$$ P(l,h) - P(l,h-h_{2}) = h(-h_{2}) + \frac{1}{2}h_{2}^{2}, $$
for all $h \geq 0$, and we get that (\ref{|q^P|_h2TlXh<C}) holds for all $l \leq 0$, all $h \geq 0$.

\noindent Finally, the inequality (\ref{|dxV|T1X1<|V|T0X0}) follows from the expression (\ref{dxV_T1X1_T0X0}) and the estimates (\ref{|q^P|_h2TlXh<C}).
\end{proof}

\noindent Let $S,m_{0,k},m_{1,k}$, $0 \leq k \leq S-1$ be positive integers.
Let $\mathcal{D}$ the linear operator from $\mathbb{C}[[\xi,\xi^{-1},x]]$ into $\mathbb{C}[[\xi,\xi^{-1},x]]$ defined by
$$ \mathcal{D}(V(\xi,x)) := \partial_{x}^{S}V(\xi,x) - \sum_{k=0}^{S-1} a_{k}(x)
(\partial_{x}^{k}V)(q^{m_{0,k}}\xi,x/q^{m_{1,k}}),$$
for all $V \in \mathbb{C}[[\xi,\xi^{-1},x]]$, where $a_{k}(x) = \sum_{s \in I_k} a_{ks}x^{s} \in \mathbb{C}[x]$, with $I_{k}$ be a finite subset
of $\mathbb{N}$, for $0 \leq k \leq S-1$.\medskip

\noindent We make the following hypothesis.\medskip

\noindent {\bf Assumption (A)} For all $0 \leq k \leq S-1$, for all $s \in I_{k}$, we have
$$ s + S - k \geq  2m_{0,k} \ \ , \ \ m_{1,k} \ge  s + S - k. $$
We consider the operator $\mathcal{A}$ from $\mathbb{C}[[\xi,\xi^{-1},x]]$ into $\mathbb{C}[[\xi,\xi^{-1},x]]$ defined by
$$ \mathcal{A}(V(\xi,x)) = V(\xi,x) - \mathcal{D}( \partial_{x}^{-S}V(\xi,x) )
 = \sum_{k=0}^{S-1} a_{k}(x)(\partial_{x}^{k-S}V)(q^{m_{0,k}}\xi,x/q^{m_{1,k}}) $$
for all $V \in \mathbb{C}[[\xi,\xi^{-1},x]]$.\medskip

\noindent From Lemma 1, we deduce the following
\begin{lemma} Let $T>0$. Then, there exists $X>0$ such that $\mathcal{A}$ is a linear bounded operator from
$(\mathbb{E}_{q,(T,X)}, ||\cdot||_{(T,X)})$ into itself. Moreover, we have that
$$ ||\mathcal{A}( V(\xi,x) )||_{(T,X)} \leq \frac{1}{2}||V(\xi,x)||_{(T,X)},$$
for all $V \in \mathbb{E}_{q,(T,X)}$.
\end{lemma}
\noindent From Lemma 4, we deduce the next
\begin{corol} Let $T>0$. Then, there exists $X>0$ such that $\mathcal{D} \circ \partial_{x}^{-S}$ is an invertible linear
operator from $(\mathbb{E}_{q,(T,X)}, ||.||_{(T,X)})$ into itself. In particular, there exists $C>0$ such that
$$
|| \mathcal{D} ( \partial_{x}^{-S}b(\xi,x) ) ||_{(T,X)} \leq C||b(\xi,x)||_{(T,X)}
$$
for all $b(\xi,x) \in \mathbb{E}_{q,(T,X)}$.
\end{corol}

\begin{defin}
Let $T_{0,j}>0$, $0 \leq j \leq S-1$ be real numbers. We say that $(T_{0,j})_{0 \leq j \leq S-1}$
satisfies {\bf the assumption (B)} if\\
1) The set
$$ \mathbb{T} := \cap_{0 \leq k \leq S-1}\cap_{k \leq j \leq S-1, s \in I_{k}}
[ q^{-m_{0,k}}T_{0,j} , T_{0,j}q^{(\frac{s+j-k - 2m_{0,k}}{2})} ] $$
is not empty,\\
2) There exists $T_{1} \in \mathbb{T}$ such that
$$ \mathbb{T}_{1} := [T_{1},T_{1}q^{S/2}] \bigcap \cap_{j=0}^{S-1} [T_{0,j},T_{0,j}q^{j/2}] $$
is not empty.
\end{defin}

\noindent {\bf Example:} Let $S \geq 1$. For all $0 \leq k \leq S-1$, let $I_{k} \subset \mathbb{N}$ such that
$s \in I_{k}$ implies $s \geq 2m_{0,k}+k$. For all $0 \leq j \leq S-1$, we put
$T_{0,j} = T_{0} > 0$. Then, we have that $T_{0} \in \mathbb{T}$ and that $\mathbb{T}_{1} = \{T_{0} \}$. So that $(T_{0,j})_{0 \leq j \leq S-1}$ satisfies the assumption {\bf (B)}.

\begin{theo}\label{teo265} Let $S \geq 1$ be an integer. For all $0 \leq k \leq S-1$, let $m_{0,k},m_{1,k}$ be positive integers and
$a_{k}(x) = \sum_{s \in I_k} a_{ks}x^{s} \in \mathbb{C}[x]$. We make the hypothesis that the assumption {\bf (A)} holds.\medskip

\noindent We consider the following functional equation
\begin{equation}
\partial_{x}^{S}V(\xi,x) = \sum_{k=0}^{S-1} a_{k}(x)(\partial_{x}^{k}V)(q^{m_{0,k}}\xi,x/q^{m_{1,k}}) \label{CP1}
\end{equation}
with initial conditions
\begin{equation}
(\partial_{x}^{j}V)(\xi,0) = \phi_{j}(\xi) \ \ , \ \ 0 \leq j \leq S-1. \label{CP1_init}
\end{equation}
We assume that $\phi_{j}(\xi) \in \mathbb{E}_{q,(T_{0,j},X_{0})}$, for $0 \leq j \leq S-1$, where $X_{0}>0$ and $(T_{0,j})_{0 \leq j \leq S-1}$ satisfy the assumption {\bf (B)}. Then, there exists $X>0$ and for $T \in \mathbb{T}_{1}$, the problem (\ref{CP1}), (\ref{CP1_init}) has a unique solution
$V(\xi,x) \in \mathbb{E}_{q,(T,X)}$. Moreover, there exists $C>0$ (depending on $S$,$q$,$a_{k}(x)$,$m_{0,k}$,$m_{1,k}$, for $0 \leq k \leq S-1$ and $X_{0}$,$T_{0,j}$, for $0 \leq j \leq S-1$) such that
$$
||V(\xi,x)||_{(T,X)} \leq C \sum_{j=0}^{S-1} ||\phi_{j}(\xi)||_{(T_{0,j},X_{0})}.
$$
\end{theo}
\begin{proof} A formal series $V(\xi,x) \in \mathbb{C}[[\xi,\xi^{-1},x]]$ which satisfies (\ref{CP1_init}) can be written in the form
$V(\xi,x)= \partial_{x}^{-S}U(\xi,x) + I(\xi,x)$ where
$$ I(\xi,x) = \sum_{j=0}^{S-1} \phi_{j}(\xi) \frac{ x^j }{ j! } $$
and $U(\xi,x) \in \mathbb{C}[[\xi,\xi^{-1},x]]$. A formal series $V(\xi,x) \in \mathbb{C}[[\xi,\xi^{-1},x]]$ is a solution of the problem
(\ref{CP1}), (\ref{CP1_init}) if and only if $U(\xi,x)$ satisfies the equation
\begin{equation}
\mathcal{D}( \partial_{x}^{-S}U(\xi,x) ) = - \mathcal{D}(I(\xi,x)). \label{eq_U}
\end{equation}
By construction, we have that
$$
 - \mathcal{D}(I(\xi,x)) = \sum_{k=0}^{S-1} \sum_{j = k}^{S-1} \sum_{s \in I_{k}} \frac{a_{ks}}{q^{m_{1,k}(j-k)}(j-k)!}
x^{s+j-k} \phi_{j}(q^{m_{0,k}}\xi)
$$
From Lemma 2 and the assumption {\bf (B)}, there exists $X_{1} >0$ such that for all $T_{1} \in \mathbb{T}$ we have that
$$ x^{s+j-k}\phi_{j}(q^{m_{0,k}}\xi) \in \mathbb{E}_{q,(T_{1},X_{1})} $$
for all $0 \leq k \leq S-1$, all $k \leq j \leq S-1$, all $s \in I_{k}$. Moreover, there exists $C_{1}>0$ (depending on $I_{k}$,$j$,$m_{0,k}$,$X_{0}$,$T_{0,j}$) such that
\begin{equation}
||x^{s+j-k}\phi_{j}(q^{m_{0,k}}\xi) ||_{(T_{1},X_{1})} \leq C_{1} ||\phi_{j}(\xi)||_{(T_{0,j},X_{0})}. \label{|xTphi_j|<|phi_j|}
\end{equation}
We deduce that $\mathcal{D}(I(\xi,x)) \in \mathbb{E}_{q,(T_{1},X_{1})}$ and from (\ref{|xTphi_j|<|phi_j|}) there exists a constant
$C_{1}'>0$ (depending on $q$,$a_{k}(x)$,$m_{0,k}$,$m_{1,k}$, for $0 \leq k \leq S-1$ and $X_{0}$,$T_{0,j}$, for $0 \leq j \leq S-1$) such that
\begin{equation}
|| \mathcal{D}(I(\xi,x)) ||_{(T_{1},X_{1})} \leq C_{1}' \sum_{j=0}^{S-1} ||\phi_{j}(\xi)||_{(T_{0,j},X_{0})}. \label{|DI|<|sumphi_j|}
\end{equation}
From Corollary 1, we deduce that the equation (\ref{eq_U}) has a unique solution $U(\xi,x) \in \mathbb{E}_{q,(T_{1},X_{1})}$. Moreover, there exists
a constant $C_{2}>0$ (depending on $q$,$a_{k}(x)$,$m_{0,k}$,$m_{1,k}$, for $0 \leq k \leq S-1$) such that
\begin{equation}
||U(\xi,x)||_{(T_{1},X_{1})} \leq C_{2}|| \mathcal{D}(I(\xi,x)) ||_{(T_{1},X_{1})}. \label{|U|<|DI|}
\end{equation}
Now, from the assumption {\bf (B)}, we choose $T_{1} \in \mathbb{T}$ in such a way that $\mathbb{T}_{1}$ is not empty. Let
$T_{2} \in \mathbb{T}_{1}$. From Lemma 3, there exists $X_{2} < X_{1}$ such that $\partial_{x}^{-S}U(\xi,x) \in
\mathbb{E}_{q,(T_{2},X_{2})}$. Moreover, there exists a constant $C_{3}>0$ (depending on $q$,$S$,$T_{1}$,$X_{1}$) such that
\begin{equation}
|| \partial_{x}^{-S}U(\xi,x) ||_{(T_{2},X_{2})} \leq C_{3} ||U(\xi,x)||_{(T_{1},X_{1})} \label{|dx-1U|<|U|}
\end{equation}
From Lemma 2, there exists $X_{3} < X_{2}$ such that $I(\xi,x) \in \mathbb{E}_{q,(T_{2},X_{3})}$. Moreover, there exists a constant
$C_{4}>0$ (depending on $S$,$q$,$T_{0,j}$,$X_{0}$, for $0 \leq j \leq S-1$) such that
\begin{equation}
|| I(\xi,x) ||_{(T_{2},X_{3})} \leq C_{4} \sum_{j=0}^{S-1} ||\phi_{j}(\xi)||_{(T_{0,j},X_{0})}. \label{|I| < |sumphi_j|}
\end{equation}
Finally, the formal series $V(\xi,x) = \partial_{x}^{-S}U(\xi,x) + I(\xi,x)$, solution of the problem (\ref{CP1}), (\ref{CP1_init}), belongs
to $\mathbb{E}_{q,(T_{2},X_{3})}$. Moreover, from the inequalities (\ref{|DI|<|sumphi_j|}), (\ref{|U|<|DI|}), (\ref{|dx-1U|<|U|}) and
(\ref{|I| < |sumphi_j|}), we get a constant $C_{5}$ (depending on $S$,$q$,$a_{k}(x)$,$m_{0,k}$,$m_{1,k}$, for $0 \leq k \leq S-1$ and $X_{0}$,$T_{0,j}$, for $0 \leq j \leq S-1$) such that
$$ ||V(\xi,x)||_{(T_{2},X_{3})} \leq C_{5} \sum_{j=0}^{S-1} ||\phi_{j}(\xi)||_{(T_{0,j},X_{0})}. $$
\end{proof}

\section{A Cauchy problem in analytic spaces of $q-$exponential growth}\label{sectqexpongrowth}

Let $S \geq 1$, $r_{1},r_{2} \geq 0$ be integers. For all $0 \leq k \leq S-1$, let $m_{0,k}$,$m_{1,k}$ be positive integers and $b_{k}(z) = \sum_{s \in I_{k}} b_{ks}z^{s}$ be a polynomial in $z$, where $I_{k}$ is a subset of $\mathbb{N}$.

\begin{lemma}\label{unique_formal_sol}
For every choice of formal series $\hat{X}_j\in\mathbb{C}[[t]]$, $0\le j\le S-1$, the Cauchy problem (\ref{qFI}), (\ref{qFI_init}) has a unique solution in the form of a formal power series $\hat{X}(t,z)=\displaystyle\sum_{h\ge 0}\hat{X}_h(t)\frac{ z^h }{ h! }$, where $\hat{X}_h\in\mathbb{C}[[t]]$ for every $h\ge  0$.
\end{lemma}

\begin{proof}
Let us put $\hat{X}_h(t)=\sum_{m\ge 0}f_{m,h}t^m$, $h\ge 0$, $f_{m,h}\in\mathbb{C}$. By substituting $\hat{X}$ in (\ref{qFI}), one can check that the left hand side turns out to be
$$
\sum_{h\ge 0}\Big(\sum_{m=0}^{r_2-1}f_{m,h+S}t^m+ \sum_{m=r_2}^{\infty}(f_{m,h+S}+(h+1)^{r_1}f_{m-r_2,h+S}q^{r_2(r_2-1)/2}q^{r_2(m-r_2)})t^m\Big)\frac{z^h}{h!},
$$
while the right hand side is
$$
\sum_{h\ge 0}\Big(\sum_{k=0}^{S-1} \sum_{h_{1}+ h_{2} = h, h_{1} \in I_{k}}
\frac{ b_{kh_{1}} t^{m_{0,k}}\hat{X}_{h_{2} + k}(q^{m_{0,k}}t) }{ h_{2}!q^{m_{1,k}h_{2}} }\Big)z^h.
$$
The values $f_{m,h}$, $m\ge 0$, $0\le h\le S-1$, are given by the initial conditions. We begin obtaining $f_{m,S}$, $m\ge 0$. We look at the coefficients of $z^0$ at both sides and impose them to be equal:
$$
\sum_{m=0}^{r_2-1}f_{m,S}t^m+ \sum_{m=r_2}^{\infty}(f_{m,S}+f_{m-r_2,S}q^{r_2(r_2-1)/2}q^{r_2(m-r_2)})t^m=
\sum_{k=0}^{S-1} \sum_{0 \in I_{k}}
b_{k0} t^{m_{0,k}}\hat{X}_{k}(q^{m_{0,k}}t).
$$
Since the right hand side is determined by the $\hat{X}_j$, $0\le j\le S-1$, we can recursively obtain the $f_{m,S}$, $m\ge 0$, by imposing the equality of the coefficients of $t^m$ in each side for every $m=0,1,2,\ldots$

\noindent We can repeat this argument in a similar way as the second index in $f_{m,h}$ increases.
\end{proof}


With the help of the $q-$Laplace transform, we reformulate our problem. Consider the Cauchy problem
\begin{equation}
((z\partial_{z} + 1)^{r_1} \tau^{r_2} + 1)
\partial_{z}^{S}\hat{W}(\tau,z) = \sum_{k=0}^{S-1} b_{k}(z)\tau^{m_{0,k}}(\partial_{z}^{k}\hat{W})(\tau,zq^{-m_{1,k}}) \label{CP2}
\end{equation}
with initial conditions
\begin{equation}
(\partial_{z}^{j}\hat{W})(\tau,0) = \hat{W}_{j}(\tau)\in\mathbb{C}[[\tau]]\ , \ \ 0 \leq j \leq S-1. \label{CP2_init}
\end{equation}

\begin{lemma}\label{equiv_Cauchy_prob}
The formal series $\hat{X}(t,z)=\displaystyle\sum_{h\ge 0}\hat{X}_h(t)\frac{ z^h }{ h! }$, where $\hat{X}_h\in\mathbb{C}[[t]]$ for every $h\ge  0$, satisfies the Cauchy problem (\ref{qFI}), (\ref{qFI_init}) if, and only if, the formal series $\hat{W}(\tau,z)=\displaystyle\sum_{h\ge 0}\hat{\mathcal{B}}_{q}\hat{X}_h(\tau)\frac{ z^h }{ h! }$ satisfies the Cauchy problem (\ref{CP2}), (\ref{CP2_init}) with $W_{j}(\tau)=\hat{\mathcal{B}}_{q}\hat{X}_j$, $0\leq j \leq S-1$.\\
Conversely, $\hat{W}(\tau,z)=\displaystyle\sum_{h\ge 0}\hat{W}_h(\tau)\frac{ z^h }{ h! }$, with  $\hat{W}_h\in\mathbb{C}[[\tau]]$ for every $h\ge  0$, satisfies the Cauchy problem (\ref{CP2}), (\ref{CP2_init}) if, and only if, the formal series $\hat{X}(t,z)=\displaystyle\sum_{h\ge 0}\hat{\mathcal{L}}_{q}\hat{W}_h(t)\frac{ z^h }{ h! }$ satisfies the Cauchy problem (\ref{qFI}), (\ref{qFI_init}) with $\hat{X}_j(t)=\hat{\mathcal{L}}_{q}\hat{W}_j(t)$ for $0\le j\le S-1$.
\end{lemma}

\begin{proof}
It suffices to insert each series in the corresponding Cauchy problem and apply (\ref{property_q_Laplace_formal}).
\end{proof}

\noindent Let $V$ be an open and bounded set in $\mathbb{C}^{\ast}$, and $q\in\mathbb{C}$ with $|q|>1$. In the following result we study the $q-$exponential growth of the coefficients of a solution to the Cauchy problem (\ref{CP2}), (\ref{CP2_init}). We will depart from initial conditions $W_j$, $0\le j\le S-1$, holomorphic in $Vq^{\mathbb{Z}}$. We make the assumptions {\bf (A)} and {\bf (B)} in the previous Section, so that we may apply Theorem~\ref{teo265}, and we also suitably choose $q$ and $V$ in order to deal with a small divisors problem.\medskip

\begin{theo}\label{theo345}
Let the assumption {\bf (A)} (of Section~\ref{sectweightedBspaces}) be fulfilled by the sets $I_{k}$ and the integers $m_{0,k}$,$m_{1,k}$, for $0 \leq k \leq S-1$.\\
1) We make the following assumptions on $q$ and on the open set $V$: $q$ is of the form $q=|q|e^{i\theta}$, with $\theta=\frac{2\pi}{br_2}$ for some $b\in\mathbb{N}$, $b\ge 1$.
If we denote
$${V}^{r_2} = \{ x^{r_2} : x \in V \}, $$
we assume that there exists $\varepsilon\in(0,\min\{\pi/b,\pi/2\})$ such that
$$
V^{r_2}\bigcap\Big(\bigcup_{l=0}^{b-1}S(-\pi+\frac{2\pi l}{b},2\varepsilon)\Big)=\emptyset,
$$
where $S(d,\varphi)$ stands for the unbounded sector in $\mathbb{C}$ with vertex at 0, bisected by direction $d$ and with opening $\varphi$.\\
2) The following assumptions on the initial conditions hold: Let $(T_{0,j})_{0 \leq j \leq S-1}$ be a sequence satisfying the assumption {\bf (B)}, there exists a constant $K_0>0$ such that
\begin{equation}
\sup_{x \in V} |W_{j}(xq^{l})| \leq K_0|q|^{\frac{1}{4}l^{2}}
(\frac{1}{T_{0,j}})^{l}\frac{1}{1+l^{2}} \ \ , \ \ \sup_{x \in V} |W_{j}(xq^{-l})| \leq
K_0(T_{0,j})^{l}\frac{1}{1+l^{2}} \label{|W_j|<K}
\end{equation}
for all $0 \leq j \leq S-1$, all $l \geq 0$.

\noindent Then, there exists a unique solution of (\ref{CP2}), (\ref{CP2_init})
$$ (\tau,z) \mapsto W(\tau,z) = \sum_{h \geq 0} W_{h}(\tau) \frac{z^h}{h!} $$
which is holomorphic on $Vq^{\mathbb{Z}} \times \mathbb{C}$.
Moreover, for all $\rho > 0$, there exist two constants $C,T>0$
(depending on $\rho$,$S$,$|q|$,$b_{k}(z)$,$m_{0,k}$,$m_{1,k}$, for $0 \leq k \leq S-1$ and
$T_{0,j}$, for $0 \leq j \leq S-1$) such that
\begin{equation}
\sup_{x \in V,z \in D(0,\rho)} |W(xq^{l},z)| \leq CK_0 |q|^{\frac{1}{2}l^{2}}(\frac{1}{T})^{l} \ \ , \ \
\sup_{x \in V,z \in D(0,\rho)} |W(xq^{-l},z)| \leq CK_0 T^{l} \label{|Wlambda|<ql2}
\end{equation}
for all $l \geq 0$ (where $K_0>0$ is defined in (\ref{|W_j|<K})).
\end{theo}

\begin{proof} From the hypothesis 1) in the statement, there exists $\delta > 0$ such that
\begin{equation}
|(h+1)^{r_1} x^{r_2}q^{r_{2}l} + 1| > \delta \label{distSV}
\end{equation}
for all $l \in \mathbb{Z}$, all $h \geq 0$, all $x \in V$. We
consider the sequence of functions $W_{h}(\tau)$, $h \geq S$,
defined as follows
\begin{equation}
\frac{ W_{h+S}(xq^{l}) }{ h! } = \sum_{k=0}^{S-1} \sum_{h_{1}+ h_{2} = h, h_{1} \in I_{k}}
\frac{ b_{kh_{1}} x^{m_{0,k}}q^{m_{0,k}l} }{ ((h+1)^{r_1} x^{r_2}q^{r_{2}l} + 1) }
\frac{ W_{h_{2} + k}(xq^{l}) }{ h_{2}!q^{m_{1,k}h_{2}} } \label{Wseq}
\end{equation}
for all $h \geq 0$, all $l \in \mathbb{Z}$, all $x \in V$. One checks that
the sequence $W_{h}(\tau)$, $h \geq 0$, of holomorphic functions on
$Vq^{\mathbb{Z}}$, satisfies the recursion (\ref{Wseq}) if and only if the formal series
$$ W(\tau,z) = \sum_{h \geq 0} W_{h}(\tau)
\frac{z^h}{h!} $$
in the $z$ variable, satisfies the problem (\ref{CP2}), (\ref{CP2_init}). From this we deduce that the solution $W$, if it exists, is unique.

\noindent According to (\ref{|W_j|<K}) and (\ref{Wseq}), we can recursively prove that the sequence $(w_{l,h})_{l \in \mathbb{Z},h \geq 0}$ defined by
\begin{equation}\label{w_lh}
w_{l,h} = \sup_{x \in V} |W_{h}(xq^{l})|,
\end{equation}
for all $l \in \mathbb{Z}$, all $h \geq 0$, consists of positive real numbers. Due to (\ref{distSV}), the sequence $(w_{l,h})_{l \in \mathbb{Z},h \geq 0}$ satisfies the following inequalities:
There exists $r>0$ (depending on $m_{0,k}$, $V$) such that
$$
\frac{ w_{l,h+S} }{ h! } \leq \sum_{k=0}^{S-1} \sum_{h_{1}+ h_{2} = h, h_{1} \in I_{k}}
\frac{ |b_{kh_{1}}| r|q|^{m_{0,k}l} }{ \delta }
\frac{ w_{l,h_{2} + k} }{ h_{2}!|q|^{m_{1,k}h_{2}} }
$$
for all $l \in \mathbb{Z}$, all $h \geq 0$.

\noindent We consider the sequence of real numbers $(v_{l,h})_{l \in \mathbb{Z},h \geq 0}$ defined by the following recursion
\begin{equation}
\frac{ v_{l,h+S} }{ h! } = \sum_{k=0}^{S-1} \sum_{h_{1}+ h_{2} = h, h_{1} \in I_{k}}
\frac{ |b_{kh_{1}}|r |q|^{m_{0,k}l} }{ \delta }
\frac{ v_{l,h_{2} + k} }{ h_{2}!|q|^{m_{1,k}h_{2}} } \label{vseq}
\end{equation}
with initial conditions $v_{l,j} = w_{l,j}$, for $0 \leq j \leq S-1$, all $l \in \mathbb{Z}$. By
construction, we have that
\begin{equation}
w_{l,h} \leq v_{l,h} \label{w<v}
\end{equation}
for all $l \in \mathbb{Z}$, all $h \geq 0$.\medskip

\noindent In the following, we put $a_{k}(x) = \sum_{s \in I_{k}} (|b_{ks}|r/\delta) x^{s}$ for $0 \leq k \leq S-1$ and we consider the formal Laurent series
$$ V(\xi,x) = \sum_{l \in \mathbb{Z},h \geq 0} v_{l,h} \xi^{l} \frac{x^{h}}{h!}. $$
From the recursion (\ref{vseq}), we get that $V(\xi,x)$ satisfies the following Cauchy problem
\begin{equation}
\partial_{x}^{S}V(\xi,x) = \sum_{k=0}^{S-1} a_{k}(x)(\partial_{x}^{k}V)(\xi |q|^{m_{0,k}},x/|q|^{m_{1,k}})
\label{CP3}
\end{equation}
with initial conditions
\begin{equation}
(\partial_{x}^{j}V)(\xi,0) = \phi_{j}(\xi) := \sum_{l \in \mathbb{Z}} w_{l,j} \xi^{l} \label{CP3_init}
\end{equation}
From the hypothesis (\ref{|W_j|<K}), we get that $\phi_{j}(\xi)$ belongs to $\mathbb{E}_{|q|,(T_{0,j},X_{0})}$, for all $X_{0}>0$. By hypothesis, the assumption {\bf (A)} holds for the sets $I_{k}$ and the numbers $m_{0,k}$,$m_{1,k}$ and the assumption
{\bf (B)} is fulfilled for the sequence $T_{0,j}$, $0 \leq j \leq S-1$. From Theorem 1, we deduce that the unique solution $V(\xi,x)$ of
the problem (\ref{CP3}), (\ref{CP3_init}) satisfies $V(\xi,x) \in \mathbb{E}_{|q|,(T,X)}$ for a real number $X>0$ and
$T \in \mathbb{T}_{1}$. Moreover,
there exists a constant $C>0$ (depending on $S$,$|q|$,$a_{k}(x)$,$m_{0,k}$,$m_{1,k}$, for $0 \leq k \leq S-1$ and $X_{0}$,$T_{0,j}$, for $0 \leq j \leq S-1$) such that
\begin{equation}
||V(\xi,x)||_{(T,X)} \leq C \sum_{j=0}^{S-1} ||\phi_{j}(\xi)||_{(T_{0,j},X_{0})}. \label{|Veps|<C|phieps|}
\end{equation}
From the inequality $P(l,h) \leq \frac{l^2}{2} - \frac{h^{2}}{4}$, for all $l \in \mathbb{Z}$, $h \geq 0$, and (\ref{|Veps|<C|phieps|})
we get that there exists a constant $C'>0$ (depending on $S$,$|q|$,$a_{k}(x)$,$m_{0,k}$,$m_{1,k}$, for $0 \leq k \leq S-1$ and $X_{0}$,$T_{0,j}$, for $0 \leq j \leq S-1$) such that
\begin{equation}
|v_{l,h}| \leq K_0C'
|q|^{\frac{l^2}{2}}|q|^{\frac{-h^{2}}{4}}h! (\frac{1}{T})^{l}(\frac{1}{X})^{h} \ \ , \ \
|v_{-l,h}| \leq K_0C'|q|^{\frac{-h^{2}}{2}} T^{l}h!(\frac{1}{X})^{h}
 \label{|vlh|<ql2q-h2}
\end{equation}
for all $l \geq 0$, all $h \geq 0$, where $K_0$ is the constant introduced in (\ref{|W_j|<K}). From the inequalities
(\ref{w<v}) and (\ref{|vlh|<ql2q-h2}), we get that
\begin{multline}
\sup_{x \in V,z \in D(0,\rho)} |W(xq^{l},z)| \leq K_0C'
|q|^{\frac{l^2}{2}} (\frac{1}{T})^{l}
(\sum_{h \geq 0} |q|^{-h^{2}/4}(\frac{\rho}{X})^{h}),\nonumber\\
\sup_{x \in V,z \in D(0,\rho)} |W(xq^{-l},z)| \leq K_0C'T^{l}(\sum_{h \geq 0} |q|^{-h^{2}/2}(\frac{\rho}{X})^{h})
\end{multline}
for all $l \geq 0$, all $\rho > 0$. So that the estimates (\ref{|Wlambda|<ql2}) hold.
\end{proof}

\noindent {\bf Remark:} Condition 1) in the previous statement could be replaced by a more general condition, namely: Let $q$ and $V$ be such that
(\ref{distSV}) is verified for some $\delta > 0$ and for all $l \in \mathbb{Z}$, all $h \geq 0$, all $x \in V$. However, we preferred to use 1) because of its easy geometrical interpretation.\medskip

\section{Second auxiliary Cauchy problem}\label{sect_second_CP}

We now suppose that the initial conditions $W_h$, $0\le h\le S-1$, of (\ref{CP2}), (\ref{CP2_init}) are holomorphic in suitably small neighbourhoods of 0. Our next aim is to obtain information on the rate of decreasing of the derivatives of the functions $W_h$, $h\ge 0$, coefficients of the solution constructed in Theorem~\ref{theo345}, near the origin. This will be done in the next Section, where we will need the second auxiliary Cauchy problem we deal with in this Section.

\begin{defin}
Let $q>1$ be given. Let us consider the space $\mathbb{H}_{(T,X)}$ of formal power series
$$
 V(\xi,x) = \sum_{l \ge 0,h \geq 0} v_{l,h} \xi^{l} \frac{x^h}{h!} \in \mathbb{C}[[\xi,x]]
$$
such that
$$|V(\xi,x)|'_{(T,X)} := \sum_{l \ge 0,h \geq 0} |v_{l,h}| T^{l} q^{h^2/2}\frac{X^h}{h!}<\infty,$$
where $T,X$ are positive real numbers.
\end{defin}
The space $(\mathbb{H}_{(T,X)}, |\cdot|'_{(T,X)})$ is a Banach algebra.

\noindent {\bf Remark:} We have a continuous inclusion $(\mathbb{H}_{(T,X')},|\cdot|'_{(T,X')})  \hookrightarrow (\mathbb{H}_{(T,X)}, |\cdot|'_{(T,X)})$
whenever $0 < X \leq X'$.\medskip

\noindent We can easily prove the following result, along the same lines as Lemma~\ref{lema1} in Section~\ref{sectweightedBspaces}.

\begin{lemma}\label{lema1tris} Let $m,s,h \geq 0$ be nonnegative integers such that $m\ge s+h$, and let $T,X>0$.
Then, there exists $C>0$ such that
$$
| x^{s}(\partial_{x}^{-h}V)(\xi, \frac{x}{q^{m}}) |'_{(T,X)} \leq C X^{s+h}
|V(\xi,x)|'_{(T,X)} \label{|xSdx-hV|<C|V|}
$$
for all $V(\xi,x) \in \mathbb{H}_{(T,X)}$.
\end{lemma}

\noindent The following is immediate from the definition of $\mathbb{H}_{(T,X)}$.
\begin{lemma}\label{lema2tris}
Let $T,X>0$. The series $R(\xi)=\sum_{l=0}^\infty 2^{l+1}\xi^l$, to be considered next, belongs to $\mathbb{H}_{(T,X)}$ if, and only if, $T<1/2$.
\end{lemma}

Let $S,m_{1,k}$, $0 \leq k \leq S-1$, be positive integers.
Let $\mathcal{F}$ be the linear operator from $\mathbb{C}[[\xi,x]]$ into $\mathbb{C}[[\xi,x]]$ defined by
$$ \mathcal{F}(V(\xi,x)) := \partial_{x}^{S}V(\xi,x) - \sum_{k=0}^{S-1} c_{k}(x)R(\xi)
(\partial_{x}^{k}V)(\xi,x/q^{m_{1,k}}),$$
for all $V \in \mathbb{C}[[\xi,x]]$, where
$$
c_k(x)=\sum_{s\in I_k}|b_{ks}|x^s,\qquad 0\le k\le S-1,
$$
and $R(\xi)$ is the one in Lemma~\ref{lema2tris}.\medskip

\noindent We consider the operator $\mathcal{B}$ from $\mathbb{C}[[\xi,x]]$ into $\mathbb{C}[[\xi,x]]$ defined by
$$ \mathcal{B}(V(\xi,x)) = V(\xi,x) - \mathcal{F}( \partial_{x}^{-S}V(\xi,x) )
 = \sum_{k=0}^{S-1} c_{k}(x)R(\xi)(\partial_{x}^{k-S}V)(\xi,x/q^{m_{1,k}}) $$
for all $V \in \mathbb{C}[[\xi,x]]$.\medskip

From now on in this Section, we will make the following\\

\noindent {\bf Assumption (A2)} For all $0 \leq k \leq S-1$, for all $s \in I_{k}$, we have
$$ m_{1,k} \ge  s + S - k. $$

From Lemmas~\ref{lema1tris} and ~\ref{lema2tris} we deduce the following result.

\begin{lemma}\label{lema3tris}
Let $T\in(0,1/2)$. Then, there exists $X>0$ such that $\mathcal{B}$ is a linear bounded operator from
$(\mathbb{H}_{(T,X)}, |\cdot|'_{(T,X)})$ into itself, and
$$ |\mathcal{B}( V(\xi,x) )|'_{(T,X)} \leq \frac{1}{2}|V(\xi,x)|'_{(T,X)},$$
for all $V \in \mathbb{H}_{(T,X)}$.
\end{lemma}
From Lemma~\ref{lema3tris}, we deduce the next

\begin{corol} Let $T\in(0,1/2)$. Then, there exists $X>0$ such that $\mathcal{F} \circ \partial_{x}^{-S}$ is an invertible linear operator from $(\mathbb{H}_{(T,X)}, |\cdot|'_{(T,X)})$ into itself. In particular, there exists $C>0$ such that
$$
| \mathcal{F} ( \partial_{x}^{-S}b(\xi,x) ) |'_{(T,X)} \leq C|b(\xi,x)|'_{(T,X)} \label{|Fdx-Sb|<|b|}
$$
for all $b(\xi,x) \in \mathbb{H}_{(T,X)}$.
\end{corol}

\begin{theo}\label{teoHCP4}
Let us consider the Cauchy problem
\begin{equation}
\partial_{x}^{S}V(\xi,x) = \sum_{k=0}^{S-1} c_{k}(x)R(\xi)(\partial_{x}^{k}V)(\xi ,x/|q|^{m_{1,k}})
\label{CP4}
\end{equation}
with initial conditions
\begin{equation}
(\partial_{x}^{j}V)(\xi,0) = \phi_{j}(\xi),\qquad 0\le j\le S-1, \label{CP4_init}
\end{equation}
and  assume that $\phi_{j}(\xi) \in \mathbb{H}_{(T_{0,j},X_{0})}$, $0 \leq j \leq S-1$, where $X_{0}>0$ and $T_{0,j}>0$ for $j=0,1,...,S-1$. Then, for every positive number $T_1<\min\{T_{0,1},\ldots,T_{0,S-1},1/2\}$, there exists $X_1>0$ such that the problem (\ref{CP4}), (\ref{CP4_init}) has a unique solution $V(\xi,x) \in \mathbb{H}_{(T_1,X_1)}$. Moreover, there exists $C>0$ (depending on $S,q,X_0$, and $c_{k}(x),m_{1,k},T_{0,k}$ for $0 \leq k \leq S-1$) such that
$$
|V(\xi,x)|'_{(T_1,X_1)} \leq C \sum_{j=0}^{S-1} |\phi_{j}(\xi)|'_{(T_{0,j},X_{0})}.
$$
\end{theo}

\begin{proof}
It heavily resembles that of Theorem~\ref{teo265},so we omit it.
\end{proof}

\section{Estimates for the derivatives of $W_j$ near the origin}\label{sectderivatorigin}

In the Cauchy problem (\ref{CP2}), (\ref{CP2_init}) we consider initial conditions $W_j$ which are holomorphic functions respectively defined in open sets containing the closed disc
$$
\overline{D}_j=\{\tau:|\tau|\le 1/( 2(j+1)^{ r_1/r_2 } )\}
$$
for $0\le j\le S-1$ (for the sake of brevity, we say that $W_j$ is holomorphic in $\overline{D}_j$). Then, Cauchy's integral formula for the derivatives allows us to obtain constants $A_j>0$ such that for every natural number $n\ge 0$ we have
$$
\max_{\tau\in \overline{D}_j}|\partial^{n}W_j(\tau)|\le A_j^n n!.
$$
So, the assumptions in the following result are not restrictive.\medskip

\begin{theo}\label{teoHCP2}
Consider the Cauchy problem (\ref{CP2}), (\ref{CP2_init}).
Suppose $W_{j}(\tau)$, $0\le j\le S-1$, are holomorphic functions in $\overline{D}_j$ such that
there exist constants $T_{0,j}>0$ and a constant $K>0$ such that
$$
\max_{\tau\in \overline{D}_{j}}|\partial^{n}W_{j}(\tau)|\le K\Big(\frac{1}{T_{0,j}}\Big)^{n}\frac{n!}{1+n^2},$$
for $n\ge0$, $j=0,1,...,S-1$.
Then there exists a formal solution of (\ref{CP2}), (\ref{CP2_init}),
$$
W(\tau,z) = \sum_{h \geq 0} W_{h}(\tau) \frac{z^h}{h!},
$$
where $W_h$ is a holomorphic function in $\overline{D}_h=\{\tau:|\tau|\leq 1/(2(h+1)^{r_1/r_2})\}$, $h\ge S$. Moreover, there exist constants $T_1,X_1>0$ such that
\begin{equation}\label{bounds_derivatives_Wh_Dhtris}
\sup_{\tau\in \overline{D}_{j}}|\partial^{n}W_{j}(\tau)|\le
C_1\Big(\frac{1}{T_1}\Big)^{n}\Big(\frac{1}{X_1}\Big)^{j}n!j!(j+1)^{r_{1}n/r_{2}}|q|^{-j^2/2},
\end{equation}
for every $n,j\ge0$, where $C_1$ is a positive
constant (depending on $S$,$q$,$b_{k}(z)$,$m_{1,k}$, for $0 \leq k
\leq S-1$ and $T_{0,j}$, for $0 \leq j \leq S-1$).
\end{theo}

\begin{proof}
We look for formal series solutions of (\ref{CP2}), (\ref{CP2_init}) of the form
$$ W(\tau,z) = \sum_{h \geq 0} W_{h}(\tau)
\frac{z^h}{h!}, $$
which leads to the equalities
\begin{equation}
\frac{ W_{h+S}(\tau) }{ h! } = \sum_{k=0}^{S-1} \sum_{h_{1}+ h_{2} = h, h_{1} \in I_{k}}
\frac{ b_{kh_{1}} \tau^{m_{0,k}} }{ ((h+1)^{r_1} \tau^{r_2} + 1) }
\frac{ W_{h_{2} + k}(\tau) }{ h_{2}!q^{m_{1,k}h_{2}} }, \label{Wseq2}
\end{equation}
for all $h\ge 0$. These equations recursively define in a unique way the sequence $\{W_h\}_{h\ge 0}$, and we easily see that $W_h$ is holomorphic in $\overline{D}_h=\{\tau:|\tau|\leq 1/(2(h+1)^{r_1/r_2})\}$, $h\ge 0$.
We aim at estimating the rate of growth of the derivatives of $W_h$ in $\overline{D}_h$.\par

Let $n_0$ be a natural number. Differentiating $n_0$ times in~(\ref{Wseq2}) we get
\begin{equation}
\frac{\partial^{n_0} W_{h+S}(\tau) }{ h! } = \sum_{k=0}^{S-1} \sum_{h_{1}+ h_{2} = h, h_{1} \in I_{k}}
b_{kh_{1}} \sum_{l_1+l_2=n_0} \frac{ n_0! }{ l_1!l_2! }\partial^{l_1}\big(\frac{  \tau^{m_{0,k}} }{ ((h+1)^{r_1} \tau^{r_2} + 1) }\big)
\frac{\partial^{l_2} W_{h_{2} + k}(\tau) }{ h_{2}!q^{m_{1,k}h_{2}} }. \label{Wseq3}
\end{equation}
It is clear that
\begin{equation}
\partial^{l_1}\big(\frac{  \tau^{m_{0,k}} }{ ((h+1)^{r_1} \tau^{r_2} + 1) }\big)=
\sum_{\lambda_1+\lambda_2=l_1,\lambda_1\le m_{0,k}}\frac{l_1!}{\lambda_1!\lambda_2!}
\frac{m_{0,k}!}{(m_{0,k}-\lambda_1)!}\tau^{m_{0,k}-\lambda_1}\partial^{\lambda_2}\big(\frac{ 1 }{ ((h+1)^{r_1} \tau^{r_2} + 1) }\big).\label{Wseq4}
\end{equation}
Following the proof of Lemma 7 in~\cite{ma1}, we get that for all $\lambda_2\ge 0$, all $\tau\in \overline{D}_h$,
\begin{equation}
\big|\partial^{\lambda_2}\big(\frac{ 1 }{ ((h+1)^{r_1} \tau^{r_2} + 1) }\big)\big|\le
\lambda_2! 2^{\lambda_2+1} (h+1)^{\frac{r_1}{r_2}\lambda_2}.\label{Wseq5}
\end{equation}
We take~(\ref{Wseq5}) into~(\ref{Wseq4}) to obtain that
\begin{eqnarray*}
\big|\partial^{l_1}\big(\frac{  \tau^{m_{0,k}} }{ ((h+1)^{r_1} \tau^{r_2} + 1) }\big)\big|&\le&
\sum_{\lambda_1+\lambda_2=l_1,\lambda_1\le m_{0,k}}\frac{l_1!}{\lambda_1!\lambda_2!}
\frac{m_{0,k}!}{(m_{0,k}-\lambda_1)!}|\tau|^{m_{0,k}-\lambda_1} \lambda_2! 2^{\lambda_2+1} (h+1)^{\frac{r_1}{r_2}\lambda_2}\\
&=&2 l_1!\sum_{\lambda_1=0}^{m_{0,k}} \frac{m_{0,k}!}{\lambda_1! (m_{0,k}-\lambda_1)!}|\tau|^{m_{0,k}-\lambda_1} 2^{l_1-\lambda_1} (h+1)^{\frac{r_1}{r_2}(l_1-\lambda_1)}\\
&=&l_1! 2^{l_1+1} (h+1)^{\frac{r_1}{r_2}l_1} \big(|\tau|+\frac{1}{2(h+1)^{r_1/r_2}}\big)^{m_{0,k}} \le l_1! 2^{l_1+1} (h+1)^{\frac{r_1}{r_2}l_1}.
\end{eqnarray*}
The previous estimates may be applied in~(\ref{Wseq3}) and they let us write
\begin{equation*}
\frac{|\partial^{n_0} W_{h+S}(\tau)| }{n_0! h! } \le \sum_{k=0}^{S-1} \sum_{h_{1}+ h_{2} = h, h_{1} \in I_{k}}
|b_{kh_{1}}| \sum_{l_1+l_2=n_0} 2^{l_1+1} (h+1)^{\frac{r_1}{r_2}l_1}
\frac{|\partial^{l_2} W_{h_{2} + k}(\tau)| }{ l_2! h_{2}! |q|^{m_{1,k}h_{2}} },
\end{equation*}
which may be rewritten as
\begin{eqnarray}
\frac{|\partial^{n_0} W_{h+S}(\tau)| }{n_0! h! (h+S+1)^{\frac{ r_1 }{ r_2 }n_0}} \le \sum_{k=0}^{S-1} \sum_{h_{1}+ h_{2} = h, h_{1} \in I_{k}}
|b_{kh_{1}}| \sum_{l_1+l_2=n_0} 2^{l_1+1}\frac{ (h+1)^{\frac{r_1}{r_2}n_0 }}{ (h+S+1)^{\frac{r_1}{r_2}n_0 } }\nonumber\\
\times \frac{ (h_2+k+1)^{\frac{r_1}{r_2}l_2 }}{ (h+1)^{\frac{r_1}{r_2}l_2 } }
\frac{|\partial^{l_2} W_{h_{2} + k}(\tau)| }{ l_2! h_{2}! (h_2+k+1)^{\frac{r_1}{r_2}l_2 } |q|^{m_{1,k}h_{2}} }.\label{Wseq7}
\end{eqnarray}
Let us put
$$
w_{n_0,h}:=\sup_{\tau\in \overline{D}_h}\frac{ |\partial^{n_0}W_h(\tau)| }{ (h+1)^{r_1n_0/r_2}},\qquad n_0\ge 0,\ h\ge 0.
$$
From~(\ref{Wseq7}) we deduce that
\begin{equation*}
\frac{ w_{n_0,h+S} }{n_0! h! } \le \sum_{k=0}^{S-1} \sum_{h_{1}+ h_{2} = h, h_{1} \in I_{k}}
|b_{kh_{1}}| \sum_{l_1+l_2=n_0} 2^{l_1+1}\frac{ w_{l_2,h_{2} + k} }{ l_2! h_{2}! |q|^{m_{1,k}h_{2}} }.
\end{equation*}
Now we define a multi-sequence $\{v_{l,h}\}_{l,h}$ by
$$
v_{l,h}=w_{l,h},\qquad l\ge 0,\ 0\le h\le S-1,
$$
and by the following recurrence relations for $n_0\ge 0$, $h\ge 0$:
\begin{equation}
\frac{ v_{n_0,h+S} }{n_0! h! } =\sum_{k=0}^{S-1} \sum_{h_{1}+ h_{2} = h, h_{1} \in I_{k}}
|b_{kh_{1}}| \sum_{l_1+l_2=n_0} 2^{l_1+1}\frac{ v_{l_2,h_{2} + k} }{ l_2! h_{2}! |q|^{m_{1,k}h_{2}} }.\label{Vseq1}
\end{equation}
It is clear that $w_{l,h}\le v_{l,h}$ for all $l\ge 0$, all $h\ge 0$. Let us consider the functions
$$
c_k(x)=\sum_{s\in I_k}|b_{ks}|x^s,\qquad 0\le k\le S-1,
$$
and
$$
R(\xi)=\sum_{l=0}^\infty 2^{l+1}\xi^l=\frac{2}{1-2\xi}.
$$
Due to the recursions  (\ref{Vseq1}), one can check that the formal power series
$\hat{V}(\xi,x)=\sum_{l,h\ge 0}\frac{v_{l,h}}{l!} \xi^l \frac{x^h}{h!}$ is a formal solution of the Cauchy problem (\ref{CP4}), (\ref{CP4_init}) with initial conditions
$$
(\partial_{x}^{j}V)(\xi,0) = \phi_{j}(\xi) := \sum_{l \ge 0} \frac{w_{l,j}}{l!} \xi^{l},\qquad 0\le j\le S-1.
$$
It is immediate to check that, for any $X_0>0$, we have $\phi_{j}(\xi)\in\mathbb{H}_{(T_{0,j},X_{0})}$ for $0\le j\le S-1$. From Theorem~\ref{teoHCP4} and the fact that $\hat{V}(\xi,x)=\sum_{n,j\ge 0}\frac{v_{n,j}}{n!} \xi^n \frac{x^j}{j!}$ is the unique formal solution of (\ref{CP4}), (\ref{CP4_init}), we can find $X_1>0$ and $T_1>0$ such that
$$|\hat{V}(\xi,x)|'_{(T_1,X_1)}\le C\sum_{j=0}^{S-1}|\phi_{j}(\xi)|'_{(T_{0,j},X_{0})}$$
for a certain $C>0$ (depending on $S$,$q$,$c_{k}(x)$,$m_{1,k}$, for $0 \leq k \leq S-1$ and $X_{0}$,$T_{0,j}$, for $0 \leq j \leq S-1$). From the last expression, one can obtain that
$$
w_{n,j}\le v_{n,j}\le C\Big(\frac{1}{T_1}\Big)^{n}\Big(\frac{1}{X_1}\Big)^{j}n!j!q^{-j^2/2}\Big[\sum_{l=0}^{S-1}|\phi_{l}(\xi)|'_{(T_{0,l},X_{0})}\Big]
\le C_1\Big(\frac{1}{T_1}\Big)^{n}\Big(\frac{1}{X_1}\Big)^{j}n!j!|q|^{-j^2/2}.
$$
We conclude by the very definition of the multi-sequence $\{w_{n,j}\}_{n,j\ge 0}$, since
$$
\max_{\tau\in \overline{D}_{j}}|\partial^{n}W_{j}(\tau)|= w_{n,j}(j+1)^{r_{1}n/r_{2}}\le C_1\Big(\frac{1}{T_1}\Big)^{n}\Big(\frac{1}{X_1}\Big)^{j}n!j!(j+1)^{r_{1}n/r_{2}}|q|^{-j^2/2},
$$
as desired.
\end{proof}

\section{Analytic solutions of the Cauchy problem with Fuchsian and irregular singularities}
\label{sect_final_solution}


Let $W_h$ be the initial data in the Cauchy problem (\ref{CP2}), (\ref{CP2_init}), and suppose they are subject to the hypotheses
of Theorem~\ref{theo345} and to the hypotheses in Theorem~\ref{teoHCP2}. Those results give us a sequence of functions $\{W_h\}_{h\ge 0}$,
holomorphic in $Vq^{\mathbb{Z}}\cup D_h$ for each $h\ge 0$, and such that the series
$$W(\tau,z) = \sum_{h \geq 0} W_{h}(\tau) \frac{z^h}{h!} $$
defines a holomorphic function on $Vq^{\mathbb{Z}} \times \mathbb{C}$ which solves the Cauchy problem.

Moreover, from (\ref{w_lh}), (\ref{w<v}) and (\ref{|vlh|<ql2q-h2}) in the proof of Theorem~\ref{theo345} we know that
\begin{equation}
\sup_{x\in V}|W_h(x q^{l})| \leq K_0C'
|q|^{\frac{l^2}{2}}|q|^{\frac{-h^{2}}{4}}h! (\frac{1}{T})^{l}(\frac{1}{X})^{h}
 \label{e855} 
\end{equation}
for all $l,h\ge 0$.

Let us choose $\lambda\in V$ and $\delta>0$.
By (\ref{e855}) we see that every $W_h$ verifies estimates as those in~(\ref{|phi|<ql2}). If we choose an integer $n(h)$ in such a way that $\lambda q^{n(h)}\in D_h$, then, according to Proposition~\ref{proposition1}, the $q-$Laplace transform of $W_h$ in the direction
$\lambda q^{n(h)} q^{\mathbb{Z}}$, which clearly equals $\lambda q^{\mathbb{Z}}$, is given by
$$\mathcal{L}_{q}^{\lambda q^{n(h)}}(W_{h})(t)=\sum_{m\in\mathbb{Z}}\frac{W_{h}(q^{m}\lambda q^{n(h)})}{\Theta(\frac{q^{m}\lambda q^{n(h)}}{t})} =\sum_{m\in\mathbb{Z}}\frac{W_{h}(q^{m}\lambda )}{\Theta(\frac{q^{m}\lambda}{t})},$$
so that it deserves to be denoted by $\mathcal{L}_{q}^{\lambda}(W_{h})(t)$. This function is well defined and holomorphic in the set $\mathcal{T}_{\lambda q^{n(h)},q,\delta,r(h)}$, which is equal to $\mathcal{T}_{\lambda,q,\delta,r(h)}$, whenever $r(h)<|\lambda q^{n(h)} q^{1/2}|T$.
We will show that these radii $r(h)$ can be taken independent of $h$, equal to $r_{0}=|\lambda q^{1/2}|T/|q|=|\lambda q^{-1/2}|T$ for every $h\ge 0$, and we will obtain precise estimates for the corresponding $q-$asymptotic expansions.

Let us assume that the function $W_h$ has the following Taylor expansion at 0,
\begin{equation}
W_h(\tau) = \sum_{n \geq 0} \frac{ f_{n,h} }{ q^{n(n-1)/2} } \tau^{n}, \label{taylor_W_h}
\end{equation}
where $f_{n,h} \in \mathbb{C}$, $n,h \geq 0$, and $\tau\in\overline{D}_h$.

\begin{prop}\label{q_Laplace_W_h}
In the situation assumed in this Section, there exist constants $B(h),D(h)>0$ (to be specified) such that
\begin{equation}\label{bounds_q_Laplace_Wh}
|\mathcal{L}_{q}^{\lambda}(W_h)(t) - \sum_{m = 0}^{n-1}f_{m,h}t^{m}|\leq D(h)B(h)^{n}|q|^{n(n-1)/2}|t|^{n}
\end{equation}
for all $n \geq 1$, for all $t \in \mathcal{T}_{\lambda,q,\delta,r_{0}}$.
\end{prop}

\begin{proof}

\noindent According to the estimates~(\ref{bounds_derivatives_Wh_Dhtris}) in Theorem~\ref{teoHCP2},
we can write
\begin{equation}\label{bound_f_nh}
\big|\frac{ f_{n,h} }{ q^{n(n-1)/2} }\big|=\big|\frac{ 1 }{ n! } \partial^n W_h(0)\big|\le C_1\Big(\frac{1}{T_1}\Big)^{n}\Big(\frac{1}{X_1}\Big)^{h}h!(h+1)^{r_{1}n/r_{2}}|q|^{-h^2/2}=
C(h)A(h)^n
\end{equation}
for every $n,h\ge 0$, where we have put, for short,
\begin{equation}\label{C(h)A(h)}
C(h)=C_1\Big(\frac{1}{X_1}\Big)^{h}h!|q|^{-h^2/2},\qquad A(h)=\frac{1}{T_1} (h+1)^{r_{1}/r_{2}},\qquad h\ge 0.
\end{equation}
For each $h\ge 0$ we define
$$
m_{h}:=\max\{m\in\mathbb{Z}:|q^m \lambda|<\frac{1}{2A(h)}\},
$$
so that
\begin{equation}\label{cota_qmlambda}
|q^{m} \lambda|A(h)<\frac{1}{2},\qquad m\le m_h.
\end{equation}
Also, we recall from Theorem~\ref{teoHCP4} that $T_1<1/2$, so
$$
\frac{1}{A(h)}<\frac{1}{2(h+1)^{\frac{r_{1}}{r_{2}}}},
$$
and we deduce that
\begin{equation}\label{qmlambdainDh}
q^m \lambda\in D_h,\qquad m\le m_h.
\end{equation}

Moreover, by the very definition of $m_h$, we have that
\begin{equation}\label{m_h+1}
m_{h}+1\ge\frac{-\log(2|\lambda|A(h))}{\log(|q|)},\qquad |q|^{m_h+1}\ge \frac{1}{2|\lambda|A(h)}.
\end{equation}
Let $K \geq 0$ be a fixed integer.
Firstly, we give estimates for
$\sum_{m > m_h} W_h(q^{m}\lambda)/\Theta(q^{m}\lambda/t)$.
Using (\ref{e855}), (\ref{prop1theta}) and (\ref{|thetalambdat|>}), we get
$$
| \frac{ W_h(q^{m}\lambda) }{ \Theta(q^{m}\lambda/t) } | \leq \frac{C'K_0h!}{K_{1}\delta}
(\frac{1}{|\lambda|})^{K} |q|^{K(K-1)/2}|t|^{K}(\frac{|t|}{T|\lambda||q|^{1/2}})^{m}\big(\frac{1}{X}\big)^h|q|^{-h^2/4}
$$
for all $m >m_h$, all $t \in \mathcal{R}_{\lambda,q,\delta}$.
For every $t\in \mathcal{T}_{\lambda,q,\delta,r_{0}}$ we have $|t|<r_{0} < T|\lambda||q|^{-1/2}$. Using (\ref{m_h+1}), we obtain that
$$
\sum_{m > m_h} (\frac{|t|}{T|\lambda||q|^{1/2}})^{m}\le \sum_{m > m_h}\big(\frac{1}{|q|}\big)^m=\frac{1}{1-|q|^{-1}}\frac{1}{|q|^{m_h+1}}\le
\frac{2|\lambda|A(h)}{1-|q|^{-1}},
$$
hence
\begin{equation}
\sum_{m > m_h} |\frac{ W_h(q^{m}\lambda) }{ \Theta(q^{m}\lambda/t) }|\leq
\frac{2|\lambda|C'K_0}{K_{1}\delta(1-|q|^{-1})}A(h)h!\big(\frac{1}{X}\big)^h|q|^{-h^2/4}
(\frac{1}{|\lambda|})^{K} |q|^{K(K-1)/2}|t|^{K} \label{|phi_m>m_h|<}
\end{equation}
for all $t \in \mathcal{T}_{\lambda,q,\delta,r_{0}}$.

\noindent In a second step, we give estimates for the sum $\sum_{m \leq m_h} W_h(q^{m}\lambda)/\Theta(q^{m}\lambda/t) -
\sum_{n = 0}^{K} f_{n,h} t^{n}$, where the $f_{n,h}$ are defined in the Taylor expansion (\ref{taylor_W_h}).
Taking into account (\ref{taylor_W_h}) and (\ref{qmlambdainDh}), we can formally write, as we did in the proof of Proposition \ref{proposition1},
\begin{multline}
\sum_{m \leq m_h} W_h(q^{m}\lambda)/\Theta(q^{m}\lambda/t) -
\sum_{n = 0}^{K} f_{n,h} t^{n} = \sum_{m \leq m_h} \frac{1}{\Theta(q^{m}\lambda/t)}\left( \sum_{n \geq K+1} \frac{f_{n,h}}{q^{n(n-1)/2}}
(q^{m}\lambda)^{n} \right) \\
- \sum_{n=0}^{K} \frac{f_{n,h}}{q^{n(n-1)/2}} \left( \sum_{m >m_h}
\frac{ (q^{m}\lambda)^{n} }{ \Theta(q^{m}\lambda/t) } \right) \label{decomp_sum_m<m_h}
\end{multline}
for all $t \in \mathbb{C}^{\ast}$.

From (\ref{decomp_sum_m<m_h}) and (\ref{bound_f_nh}), we deduce that
\begin{equation}
| \sum_{m \leq m_h} W_h(q^{m}\lambda)/\Theta(q^{m}\lambda/t) -
\sum_{n = 0}^{K} f_{n,h} t^{n} | \leq \mathcal{A}(t) + \mathcal{B}(t) \label{|phi_m<0 - polyn|<bis}
\end{equation}
where
$$ \mathcal{A}(t) = \sum_{m \leq m_h}
\frac{1}{|\Theta(q^{m}\lambda/t)|}\left( \sum_{n \geq K+1} C(h)A(h)^{n}(|q|^{m}|\lambda|)^{n} \right) $$
and
$$ \mathcal{B}(t) = \sum_{n=0}^{K} C(h)A(h)^{n} \left( \sum_{m > m_h}
\frac{ |(q^{m}\lambda)^{n}| }{ |\Theta(q^{m}\lambda/t)| } \right), $$
for all $t \in \mathbb{C}^{\ast}$.\par

We give estimates for $\mathcal{A}(t)$.
Taking into account (\ref{cota_qmlambda}), we have that
\begin{equation}
\mathcal{A}(t) \leq C(h)\sum_{m \leq m_h}\frac{1}{|\Theta(q^{m}\lambda/t)|}\frac{(A(h)|q|^m|\lambda|)^{K+1}}{1-A(h)|q|^m|\lambda|}
\le\frac{C(h)A(h)^{K+1}}{1-A(h)|q|^{m_h}|\lambda|}\sum_{m \geq -m_h} \frac{ (|q|^{-m}|\lambda|)^{K+1} }{ |\Theta(q^{-m}\lambda/t)| } \label{A1<bis}
\end{equation}
for all $t \in \mathbb{C}^{\ast}$.
From (\ref{prop2theta}), we have that
$$ |\Theta(q^{-m}\lambda/t)| \geq K_{1}\delta |q|^{-K(K-1)/2} |\frac{q^{-m}\lambda}{t}|^{K} $$
for all $m \geq -m_h$, all $t \in \mathcal{R}_{\lambda,q,\delta}$. We deduce that
\begin{equation}
\frac{ (|q|^{-m}|\lambda|)^{K+1} }{ |\Theta(q^{-m}\lambda/t)| } \leq \frac{|\lambda|}{K_{1} \delta} |q|^{K(K-1)/2}|t|^{K}
(\frac{1}{|q|})^{m} \label{A2<bis}
\end{equation}
for all $m \geq -m_h$, all $t \in \mathcal{R}_{\lambda,q,\delta}$. From (\ref{A1<bis}), (\ref{A2<bis}) and (\ref{cota_qmlambda}), we conclude that
\begin{eqnarray}
\mathcal{A}(t) &\leq &\frac{C(h)A(h)}{1-A(h)|q|^{m_h}|\lambda|}\frac{|\lambda|}{K_{1} \delta}\frac{|q|^{m_h}}{1-|q|^{-1}}A(h)^{K}|q|^{K(K-1)/2}|t|^{K}\nonumber\\
&\leq &\frac{C(h)}{K_{1} \delta(1-|q|^{-1})}A(h)^{K}|q|^{K(K-1)/2}|t|^{K} \label{calA<bis}
\end{eqnarray}
for all $t \in \mathcal{R}_{\lambda,q,\delta}$.\par

In the next step, we get estimates for $\mathcal{B}(t)$.
From (\ref{prop2theta}), we have that
$$ |\Theta(q^{m}\lambda/t)| \geq K_{1}\delta |q|^{-(K+1)K/2} |\frac{q^{m}\lambda}{t}|^{K+1} $$
for all $m > m_h$, all $t \in \mathcal{R}_{\lambda,q,\delta}$. We deduce that
$$
\frac{ |(q^{m}\lambda)^{n}| }{ |\Theta(q^{m}\lambda/t)| } \leq \frac{ |\lambda|^{n} }{ K_{1}\delta }
(\frac{1}{|\lambda|})^{K+1}|q|^{(K+1)K/2}|t|^{K+1} (\frac{1}{|q|^{K+1-n}})^{m}
$$
for all $m >m_h$, all $0 \leq n \leq K$. Then,
\begin{equation}
\sum_{m > m_h}
\frac{ |(q^{m}\lambda)^{n}| }{ |\Theta(q^{m}\lambda/t)| }\le \frac{ |\lambda|^{n} }{ K_{1}\delta }
(\frac{1}{|\lambda|})^{K+1}|q|^{(K+1)K/2}|t|^{K+1} \frac{(\frac{1}{|q|^{K+1-n}})^{m_h+1}}{1-\frac{1}{|q|^{K+1-n}}}.\label{B1<bis}
\end{equation}
It is clear that $|q|^{K+1-n}\ge |q|$ for $0\le n\le K$, hence
$$
1-\frac{1}{|q|^{K+1-n}}\ge 1-|q|^{-1},\quad 0\le n\le K.
$$
If we write
$$
(\frac{1}{|q|^{K+1-n}})^{m_h+1}=\big(\frac{1}{|q|^{m_h+1}}\big)^{K+1}(|q|^{m_h+1})^{n},
$$
from (\ref{B1<bis}) we deduce that
\begin{align*}
\mathcal{B}(t)&\le \frac{C(h)}{K_1\delta(1-|q|^{-1})}(\frac{1}{|\lambda||q|^{m_h+1}})^{K+1}|q|^{(K+1)K/2}|t|^{K+1}\sum_{n=0}^{K} (A(h)|\lambda||q|^{m_h+1})^{n}\\
&\le \frac{C(h)}{K_1\delta(1-|q|^{-1})}(\frac{1}{|\lambda||q|^{m_h+1}})^{K+1}|q|^{(K+1)K/2}|t|^{K+1}\sum_{n=0}^{K} (2A(h)|\lambda||q|^{m_h+1})^{n}.
\end{align*}
By (\ref{m_h+1}) we know that $2A(h)|\lambda||q|^{m_h+1}\ge 1$. For every real number $x\ge 1$ we have
$$
\sum_{n=0}^K x^n\le \sum_{n=0}^K {K\choose n}x^nx^{K-n}=(2x)^K,
$$
and we deduce that $\sum_{n=0}^{K} (2A(h)|\lambda||q|^{m_h+1})^{n}\le (4A(h)|\lambda||q|^{m_h+1})^{K}$. Also, we have $|t|<\!|\lambda q^{-1/2}|T$ whenever $t \in \mathcal{T}_{\lambda,q,\delta,r_{0}}$, and
$(K+1)K/2=K+K(K-1)/2$. Gathering all these facts and using (\ref{m_h+1}), we deduce that
\begin{eqnarray}
\mathcal{B}(t)&\le &\frac{ |t| C(h) }{ K_1 \delta (1-|q|^{-1}) |\lambda| |q|^{m_h+1} }(4A(h)|q|)^{K}|q|^{K(K-1)/2}|t|^{K}\nonumber\\
&\le& \frac{ 2 |\lambda| T A(h) C(h) }{ K_1 \delta |q|^{1/2} (1-|q|^{-1}) } ( 4 A(h) |q| )^{K} |q|^{K(K-1)/2}|t|^{K}\label{calB<bis}
\end{eqnarray}
for all $t \in \mathcal{T}_{\lambda,q,\delta,r_{0}}$.\medskip

\noindent Finally, using the estimates
\begin{multline*}
| \sum_{m \in \mathbb{Z}} W_h(q^{m}\lambda)/\Theta(q^{m}\lambda/t) -\sum_{n = 0}^{K} f_{n,h} t^{n} |
\leq | \sum_{m \leq m_h} W_h(q^{m}\lambda)/\Theta(q^{m}\lambda/t) -\sum_{n = 0}^{K} f_{n,h} t^{n} | \\
+ | \sum_{m > m_h} W_h(q^{m}\lambda)/\Theta(q^{m}\lambda/t) |
\end{multline*}
we deduce from (\ref{|phi_m>m_h|<}), (\ref{|phi_m<0 - polyn|<bis}), (\ref{calA<bis}), (\ref{calB<bis}) that
\begin{equation}\label{bound_1_q_Laplace_Wh}
|\mathcal{L}_{q}^{\lambda}(W_h)(t) - \sum_{n = 0}^{K}f_{n,h}t^{n}|\leq D_1(h)B_1(h)^{K}|q|^{K(K-1)/2}|t|^{K}
\end{equation}
for all $K \geq 0$, for all $t \in \mathcal{T}_{\lambda,q,\delta,r_{0}}$, with
\begin{equation}\label{constants_1_Bh_Dh}
B_1(h)=B_1(h+1)^{r_1/r_2},\qquad D_1(h)=B_2 (h+1)^{r_1/r_2} h! B_3^h |q|^{-h^2/4},
\end{equation}
where $B_1, B_2$ and $B_3$ are positive constants that do not depend on $h$.
In order to conclude, it suffices to write, for $K\ge 1$,
$$
|\mathcal{L}_{q}^{\lambda}(W_h)(t) - \sum_{n = 0}^{K-1}f_{n,h}t^{n}|\leq
|\mathcal{L}_{q}^{\lambda}(W_h)(t) - \sum_{n = 0}^{K}f_{n,h}t^{n}|+|f_{K,h}t^K|,
$$
and take into account (\ref{bound_1_q_Laplace_Wh}) and (\ref{bound_f_nh}). According to the expressions
(\ref{C(h)A(h)}) and (\ref{constants_1_Bh_Dh}), one obtains the estimates (\ref{bounds_q_Laplace_Wh})
with
\begin{equation}\label{constants_Bh_Dh}
B(h)=A_1(h+1)^{r_1/r_2},\qquad D(h)=A_2 (h+1)^{r_1/r_2} h! A_3^h |q|^{-h^2/4},
\end{equation}
where $A_1, A_2$ and $A_3$ are again positive constants that do not depend on $h$.
\end{proof}

\noindent We are ready to obtain our main result.

\begin{theo}\label{teo_solution}
Suppose $\hat{X}_j(t)=\sum_{m\ge 0}f_{m,j}t^m\in\mathbb{C}[[t]]$, $0\le j\le S-1$, are given initial conditions for the Cauchy problem (\ref{qFI}), (\ref{qFI_init}), and let
$$
\hat{X}(t,z)=\sum_{h\ge 0}\hat{X}_h(t)\frac{z^h}{h!}=\sum_{h\ge 0}\sum_{m\ge 0}f_{m,h}t^m\frac{z^h}{h!}
$$
be the only formal series solution of the problem (see Lemma~\ref{unique_formal_sol}). We suppose that the series $\hat{X}_j(t)$, $0\le j\le S-1$, are $q-$Gevrey of order 1, and that their formal $q-$Borel transforms of order~1, $W_j(\tau)=\hat{\mathcal{B}}_{q}\hat{X}_j(\tau)$, which are holomorphic functions around 0, indeed satisfy the assumptions of Theorems~\ref{theo345} and~\ref{teoHCP2}. We also assume that the rest of hypotheses of Theorem~\ref{theo345} are satisfied. Let
$$
W(\tau,z)=\sum_{h \geq 0} W_{h}(\tau) \frac{z^h}{h!}
$$
be the solution of the Cauchy problem (\ref{CP2}), (\ref{CP2_init}), corresponding to the initial conditions $W_j$, $0\le j\le S-1$. Then, we have that:\\
\noindent 1) The function $X(t,z)=\displaystyle\sum_{h\ge 0}\mathcal{L}_{q}^{\lambda}(W_h)(t)\frac{z^h}{h!}$ is holomorphic in $\mathcal{T}_{\lambda,q,\delta,r_{0}}\times \mathbb{C}$.\\
\noindent 2) The function $X(t,z)$ solves the Cauchy problem (\ref{qFI}), (\ref{qFI_init}).\\
\noindent 3) If $r_1\ge 1$, given $R>0$ there exist constants $\tilde{C}>0$, $\tilde{D}>0$ such that for every $n\in\mathbb{N}$, $n\ge 1$, one has
\begin{equation}\label{bounds_solution}
\Big|X(t,z)-\sum_{h\ge 0}\sum_{m=0}^{n-1}f_{m,h}t^m\frac{z^h}{h!}\Big|\le \tilde{C}\tilde{D}^n\Gamma(\frac{r_1}{r_2}(n+1))|q|^{n(n-1)/2}|t|^n
\end{equation}
for every $t\in\mathcal{T}_{\lambda,q,\delta,r_{0}}$, $z\in D(0,R)$.\\
\noindent If $r_1=0$, given $R>0$ there exist constants $\tilde{C}>0$, $\tilde{D}>0$ such that for every $n\in\mathbb{N}$, $n\ge 1$, one has
\begin{equation}\label{bounds_solution_2}
\Big|X(t,z)-\sum_{h\ge 0}\sum_{m=0}^{n-1}f_{m,h}t^m\frac{z^h}{h!}\Big|\le \tilde{C}\tilde{D}^n|q|^{n(n-1)/2}|t|^n
\end{equation}
for every $t\in\mathcal{T}_{\lambda,q,\delta,r_{0}}$, $z\in D(0,R)$.
\end{theo}

\noindent {\bf Remark:} Due to the estimates (\ref{bounds_solution}) and (\ref{bounds_solution_2}), we may say that the function $X(t,z)$ admits the series $\displaystyle\sum_{h\ge 0}\sum_{m\ge 0}f_{m,h}t^m\frac{z^h}{h!}$ as
$q-$asymptotic expansion of order 1 in $t$, uniformly for $z$ in the compact subsets of $\mathbb{C}$. It may be noted that, because of the small divisors problem we have dealt with, a new factor appears in the estimates, in terms of the Eulerian Gamma function. The value $r_1/r_2$ may be thought of as a sub-order, or a second-level order, in the asymptotic expansion.\medskip

\begin{proof}
1) In view of (\ref{bounds_q_Laplace_Wh}), for $n=1$, and (\ref{constants_Bh_Dh}) we have that
$$
|\mathcal{L}_{q}^{\lambda}(W_h)(t) - f_{0,h}|\leq D(h)B(h)|t|\le
A_1(h+1)^{2r_1/r_2}A_2 h! A_3^h |q|^{-h^2/4}r_0
$$
for every $h\ge 0$, every $t\in\mathcal{T}_{\lambda,q,\delta,r_{0}}$. On the other hand, by (\ref{bound_f_nh}) we have that
$$
|f_{0,h}|\le C_1\Big(\frac{1}{X_1}\Big)^{h}h! |q|^{-h^2/2}
$$
for every $h\ge 0$. So, we conclude that there exist $A_4,A_5>0$ such that
$$
|\mathcal{L}_{q}^{\lambda}(W_h)(t)|\leq A_4 A_5^h h! |q|^{ -h^2/4}
$$
for every $h\ge 0$, every $t\in\mathcal{T}_{\lambda,q,\delta,r_{0}}$. Then, for $z\in D(0,R)$ we have
$$
|\sum_{h\ge 0}\mathcal{L}_{q}^{\lambda}(W_h)(t)\frac{z^h}{h!}|\le \sum_{h\ge 0} A_4(A_5R)^h |q|^{-h^2/4}<\infty,
$$
so that the series converges and the function it defines is holomorphic in $\mathcal{T}_{\lambda,q,\delta,r_{0}}\times \mathbb{C}$.\\
\noindent 2) Since the series $\sum_{h \geq 0} W_{h}(\tau) \displaystyle\frac{z^h}{h!}$ is a
solution of (\ref{CP2}), (\ref{CP2_init}), one can guarantee that $X(t,z)$ is a solution of the Cauchy problem (\ref{qFI}), (\ref{qFI_init}) by Proposition \ref{qLaplacetauphi}.\\
\noindent 3) For every $n\ge 1$, every $(t,z)\in\mathcal{T}_{\lambda,q,\delta,r_{0}}\times D(0,R)$, the sum
$$
\sum_{h\ge 0}\sum_{m=0}^{n-1}f_{m,h}t^m\frac{z^h}{h!}
$$
is convergent, as we see from~(\ref{bound_f_nh}). One may take into account~(\ref{bounds_q_Laplace_Wh}) and~(\ref{constants_Bh_Dh}) and write
\begin{eqnarray}
\Big|X(t,z)-\sum_{h\ge 0}\sum_{m=0}^{n-1}f_{m,h}t^m\frac{z^h}{h!}\Big|&\le&
\sum_{h\ge 0} |\mathcal{L}_{q}^{\lambda}(W_h)(t)-\sum_{m=0}^{n-1}f_{m,h}t^m|\frac{R^h}{h!}\nonumber\\
&\le& A_2 A_1^n |q|^{n(n-1)/2}|t|^{n} \sum_{h\ge 0} (h+1)^{r_1(n+1)/r_2} (A_3R)^h |q|^{-h^2/4}\nonumber\\
&=&\frac{A_2}{A_3R} A_1^n |q|^{n(n-1)/2}|t|^{n} \sum_{h\ge 1} h^{r_1(n+1)/r_2} (A_3R)^h |q|^{-(h-1)^2/4}.\label{sum_Hankel}
\end{eqnarray}
In case $r_1=0$, the conclusion easily follows, since the last sum is convergent and independent of $n$. In case $r_1\ge 1$, we follow an idea of B. Braaksma and L. Stolovitch~\cite{BraaStol}. Let $\varepsilon>0$, and let $\gamma$ be a contour that goes from $\infty e^{-i\pi}$ to $-\varepsilon$ along the negative real axis, then it turns once around 0 in the positive sense, and it goes from $-\varepsilon$ to $\infty e^{i\pi}$ again along the negative real axis. For
\begin{equation}\label{mu}
\mu=\frac{r_1(n+1)}{r_2}>0,
\end{equation}
Hankel's formula allows us to write
$$
\frac{h^\mu}{\Gamma(\mu+1)}=\frac{1}{2\pi i}\int_{\gamma}e^{hs}s^{-\mu-1}\,ds,
$$
so that the sum in~(\ref{sum_Hankel}) may be written as
\begin{multline}
\frac{\Gamma(\mu+1)}{2\pi i} \sum_{h\ge 1} (A_3R)^h |q|^{-(h-1)^2/4}\int_{\gamma}e^{hs}s^{-\mu-1}\,ds\\
=\frac{\Gamma(\mu+1)}{2\pi i} \sum_{h\ge 1}\int_{\gamma}s^{-\mu-1}|q|^{-(h-1)^2/4} (A_3Re^s)^h \,ds.\label{sum_integral_Hankel}
\end{multline}
We consider now the entire function
$$
F(z)=\sum_{h\ge 1} |q|^{-(h-1)^2/4}z^h,\qquad z\in\mathbb{C}.
$$
The series converges uniformly in every closed disc. Observe that as $s$ runs over $\gamma$, its real part remains bounded above, and the same is valid for the modulus of
$A_3Re^s$. So, we may write
$$
F(A_3Re^s)=\sum_{h\ge 1} |q|^{-(h-1)^2/4}(A_3Re^s)^h
$$
uniformly in $\gamma$, and the dominated convergence theorem ensures that
\begin{multline}\label{integral_F}
\sum_{h\ge 1}\int_{\gamma}s^{-\mu-1} |q|^{-(h-1)^2/4}(A_3Re^s)^h \,ds=
\int_{\gamma}s^{-\mu-1} \sum_{h\ge 1}|q|^{-(h-1)^2/4}(A_3Re^s)^h \,ds\\=
\int_{\gamma}s^{-\mu-1} F(A_3Re^s)\,ds.
\end{multline}
Moreover, $F(A_3Re^s)$ remains bounded as $s$ runs over $\gamma$, say by $M>0$,
and it is easy to obtain, estimating on each of the three parts of $\gamma$, that
\begin{equation}\label{bound_integral_F}
|\int_{\gamma}s^{-\mu-1} F(A_3Re^s)\,ds|\le 2\frac{M}{\mu \varepsilon^{\mu}}+\frac{2\pi M}{\varepsilon^{\mu}}\le\frac{\tilde{M}^{\mu}}{\mu \varepsilon^{\mu}},
\end{equation}
where $\tilde{M}>0$ is some suitable constant independent of $h$.
Gathering (\ref{sum_Hankel}), (\ref{sum_integral_Hankel}), (\ref{integral_F}) and (\ref{bound_integral_F}), we see that
$$
\Big|X(t,z)-\sum_{h\ge 0}\sum_{m=0}^{n-1}f_{m,h}t^m\frac{z^h}{h!}\Big|\le
\frac{A_2}{A_3R} A_1^n |q|^{n(n-1)/2}|t|^{n} \frac{\Gamma(\mu+1)}{2\pi i}
\frac{\tilde{M}^{\mu}}{\mu \varepsilon^{\mu}}.
$$
It suffices to recall that $\Gamma(\mu+1)=\mu\Gamma(\mu)$ and the definition of $\mu$, (\ref{mu}), in order to conclude.
\end{proof}

\noindent{\bf Remark:} All the results in this work are valid for any $r_1\ge 0$, but the case $r_1=0$, as it may be seen in the last Theorem, deserves some attention, since the Fuchsian singularity at $z=0$ does not appear any more. The most important consequence of this fact is the disappearance of the small divisors phenomenon we had in general.\\
\noindent Moreover, the condition 1) in Theorem~\ref{theo345}, concerning the argument of $q$ and the set $V$, can be relaxed. Indeed, the estimates~(\ref{distSV}) hold if one assumes that there exists $\delta>0$ such that
\begin{equation}\label{condition_V_q2}
\mathrm{dist}(V^{r_2}q^{r_2\mathbb{Z}},\{-1\})>\delta,
\end{equation}
where $\mathrm{dist}$ is the Euclidean distance between two sets in $\mathbb{C}$. For example,
suppose $V$ is such that there exist $R_1,R_2$ with
$$
0<R_1\le |x^{r_2}|\le R_2
$$
for all $x\in V$, and suppose that $R_2<|q|R_1$ and there exists $j\in\mathbb{Z}$ such that
$$
|q|^{r_2j}R_2<1<|q|^{r_2(j+1)}R_1.
$$
Then, one can easily check that the condition~(\ref{condition_V_q2}) holds.\\
\noindent In Theorem~\ref{teoHCP2} all the functions $W_h$ are holomorphic in a common disc, say  $\overline{D}$, and there exist constants $T_1,X_1>0$ such that
$$
\sup_{\tau\in \overline{D}}|\partial^{n}W_{j}(\tau)|\le
C_1\Big(\frac{1}{T_1}\Big)^{n}\Big(\frac{1}{X_1}\Big)^{j}n!j!|q|^{-j^2/2}
$$
for every $n,j\ge 0$. The proof of Proposition~\ref{q_Laplace_W_h} admits some simplification, and one obtains that
$$
|\mathcal{L}_{q}^{\lambda}(W_h)(t) - \sum_{m = 0}^{n-1}f_{m,h}t^{m}|\leq
A_2 h! A_3^h |q|^{-h^2/4}A_1^{n}|q|^{n(n-1)/2}|t|^{n},
$$
for every $h\ge 0$, $n\ge 1$. Finally, no sub-order appears in the $q-$asymptotic expansion of the solution $X(t,z)$.

%

\end{document}